\documentclass[10pt, reqno, a4paper]{amsart}

\usepackage{amsthm,amsmath,amsfonts,graphicx}
\usepackage{a4wide}
\usepackage{bm}
\usepackage{thmtools}
\usepackage{xcolor}
\usepackage{enumitem}

\theoremstyle{plain}
\declaretheorem[name=Proposition, numberwithin=section]{Proposition}
\declaretheorem[name=Theorem, sibling=Proposition]{Theorem}
\declaretheorem[name=Lemma, sibling=Proposition]{Lemma}

\theoremstyle{definition}
\declaretheorem[name=Remark, sibling=Proposition]{Remark}

\usepackage[colorlinks, linkcolor=blue]{hyperref}

\renewcommand{\P}{\ensuremath{\mathbb{P}\,}}
\newcommand{\R}{\ensuremath{\mathbb{R}}}
\newcommand{\N}{\ensuremath{\mathbb{N}}}
\newcommand{\F}{\ensuremath{\mathcal{F}}}
\renewcommand{\d}{\ensuremath{\,{\rm d}}}
\newcommand{\e}{\ensuremath{{\rm e}}}

\newcommand\numberthis{\addtocounter{equation}{1}\tag{\theequation}}

\makeatletter
\@namedef{subjclassname@2020}{\textup{2020} Mathematics Subject Classification}
\makeatother

\begin{document}

%\begin{frontmatter}

% "Title of the Paper"
\title[Parameter estimation for SDEs with Rosenblatt processes]{Parameter estimation and singularity of laws on the path space for SDEs driven by Rosenblatt processes}

\author{Petr \v Coupek}
\address{Charles University, Faculty of Mathematics and Physics, Sokolovsk\' a 83, Prague 8, 186 75, Czech Republic. }
\email{coupek@karlin.mff.cuni.cz}

\author{Pavel K\v r\'i\v z}
\address{Charles University, Faculty of Mathematics and Physics, Sokolovsk\' a 83, Prague 8, 186 75, Czech Republic.}
\email{kriz@karlin.mff.cuni.cz}

\author{Bohdan Maslowski}
\address{Charles University, Faculty of Mathematics and Physics, Sokolovsk\' a 83, Prague 8, 186 75, Czech Republic.}
\email{maslow@karlin.mff.cuni.cz}

\thanks{This research was supported by the Czech Science Foundation
project No.\ 22-12790S}


% 60G22 Fractional processes, including fractional Brownian motion
% 62M09 Non-Markovian processes: estimation

\begin{abstract}
In this paper, we study parameter identification for solutions to (possibly non-linear) SDEs driven by additive Rosenblatt process and singularity of the induced laws on the path space. We propose a joint estimator for the drift parameter, diffusion intensity, and Hurst index that can be computed from discrete-time observations with a bounded time horizon and we prove its strong consistency (as well as the speed of convergence) under in-fill asymptotics with a fixed time horizon. As a consequence of this strong consistency, singularity of measures generated by the solutions with different drifts is shown. This results in the invalidity of a Girsanov-type theorem for Rosenblatt processes.
\end{abstract}

\maketitle

\section{Introduction}

Models with self-similar processes that exhibit long-range dependence experience noticeable interest among researchers. The need for such models emerge in many fields (such as network traffic, mathematical finance, hydrology, image processing, etc.) and we refer, for example, to \cite{Tud13} or the more recent monographs \cite{PipiTaqqu17,Sam16} for a thorough exposition.

The most popular family of self-similar processes with long-range dependence is probably the family of fractional Brownian motions (FBMs) with Hurst parameter $H\in (1/2,1)$ as these processes are Gaussian and have stationary increments. We refer to e.g., \cite{DecrUstu99,MvN68} for its definition and basic properties and to \cite{BiaHuOksZha08} and to the references therein for a more detailed account.

If, however, Gaussianity is not a reasonable assumption, other choices for the driving process must be considered. A more general class of self-similar processes with long-range dependence consists of the so-called Hermite processes; see, e.g., \cite{Tud13}, and the closest family to that of FBMs within this class is the family of Rosenblatt processes (FBMs are Hermite processes of order $1$ while Rosenblatt processes are Hermite processes of order $2$). These processes can be defined as iterated Wiener integrals (w.r.\ to the standard Wiener process) of a deterministic kernel (which in particular means that they are not Gaussian unless they are of order one) and they arise naturally as limits of suitably normalized sums of strongly dependent random variables. We refer to the survey article \cite{Taqq11} for a very accessible exposition of their construction and basic properties, to \cite{Tud08} for their thorough analysis, and to \cite{AbrPip06,Alb98,CouOnd22,GarTorTud12,KerNouSakVii21,Pip04} for some of their finer properties. We note as well that Rosenblatt processes are also considered as the driving noise for stochastic (partial) differential equations (S(P)DEs) whose analysis has undergone vital development in recent years; see, e.g., \cite{BonTud11,CouMasOnd18,CouMasOnd22,SlaTud19,SlaTud19b,SlaTud19c}.

Should models with self-similar processes with long-range dependence be applicable in practice, their efficient calibration is a must. The literature on parameter estimation for SDEs driven by FBMs is already quite extensive and a comprehensive reading that includes numerous references can be found for example in the book \cite{Kuto04} for diffusion processes, and in the monographs \cite{KubiliusEtAl2017} or \cite{PrakasaRao10} for fractional diffusion processes. 
A notable approach - that of maximum likelihood estimation - is based on a Girsanov-type theorem for FBMs. This theorem, established in \cite{DecrUstu99}, implies that the probability laws of the solutions to SDEs driven by FBMs corresponding to different values of the drift parameter on the space of continuous trajectories are equivalent. Based on this theorem, one can, for example, obtain a consistent drift estimator if the time horizon is infinite; see, e.g., papers \cite{HuNualartZhou19} for the linear case and \cite{TudorViens07} for the non-linear case. In contrast, consistent estimation of the diffusion coefficient and the Hurst parameter of the driving noise is possible even under in-fill asymptotics with a fixed time horizon; see \cite{BerzinEtAl2014}.

While parameter estimation for models with FBMs has been much investigated, parameter estimation for models with Rosenblatt (or other higher-order Hermite) processes is still in its infancy and the literature on this topic is rather sparse. An estimator of the Hurst parameter $H$ for the standardized Rosenblatt process $\{Z^H(t)\}_{t \in [0,1]}$ in a high-frequency regime is presented in \cite{TudorViens09}. The estimator considered therein is based on $2$-variations. Estimators for $H$ that are based on $2$-variations for longer filters are further studied in \cite{ChronopoulouEtAl09}.
A generalization of the results from \cite{TudorViens09} to higher-order Hermite processes can be found in \cite{ChronopoulouEtAl11}. Estimation of $H$ for a rescaled Rosenblatt process $\sigma Z^H$ that uses wavelet coefficients is studied in \cite{BarTud10}. As far as parameter estimation in SDEs with Rosenblatt noise is concerned, the 2-variation approach from \cite{TudorViens09} is used to estimate $H$ of the driving Hermite process in the Hermite Ornstein-Uhlenbeck model (solution to the Langevin equation) in \cite{AssaadTudor20}. However, in that paper, the presence of the unknown scaling parameter for the driving process significantly reduces the speed of convergence of the estimator to the true value of the parameter. In the same paper, an estimator of the scale parameter $\sigma$ based on the generalized variation with known $H$ is also presented and studied. There are also estimators of the drift parameter. In \cite{NourdinTran19}, the authors study a moment-type estimator of the drift parameter of a Hermite Ornstein-Uhlenbeck process that is based on a continuously (in time) observed trajectory and they prove strong consistency of the estimator with the time horizon tending to infinity. A consistent maximum-likelihood type estimator of the drift of a drifted Rosenblatt process is also given in \cite{BerTorTud11}. The estimator is not constructed via a Girsanov-type theorem, that seems not to have been investigated so far, but by approximation of the Rosenblatt process by a two-dimensional disturbed random walk, that converges to the Rosenblatt process in the Skorokhod topology, instead.

In this paper, we focus on parameter estimation in non-linear SDEs with additive Rosenblatt noise and we address the closely related question on equivalence/singularity of the measures on the path space corresponding to different drift parameters. Our aim here is to provide a joint estimation procedure for the three main parameters (the drift parameter, diffusion intensity, and the self-similarity index) which is strongly consistent in a high-frequency regime with in-fill asymptotics. More precisely, we consider the non-linear SDE with additive noise given by
	\begin{equation}
		\label{SDE}
		\d X(t) =  \lambda f(X(t)) \d t + \sigma \d Z^H(t),\quad X(0)=X_0,
	\end{equation}
where
\begin{itemize}
\itemsep0em
\item $\lambda \in \R$ is an unknown drift parameter,
\item $f: \R \to \R$ is locally Lipschitz and satisfies a certain Lyapunov-type condition,
\item $\sigma > 0$ is an unknown scaling (noise intensity) parameter,
\item $\{Z^H(t)\}_{t \in [0,1]}$ is a Rosenblatt process with unknown Hurst parameter $H \in (1/2,1)$, 
\item $X_0\in\mathbb{R}$ is a known initial condition.
\end{itemize}
In order to estimate the values of $\lambda, \sigma$, and $H$, we observe a single trajectory of the solution $X$ sampled at discrete time instants, i.e. our data is a finite sequence $\{X(i/N): i = 0,1,\ldots,N\}$ for some $N\in\N$. We consider discrete observations on the fixed time-window $[0,1]$ with decreasing mesh size ($1/N \to 0$). This in-fill asymptotics is appropriate for high-frequency data.

After some necessary preliminaries in section \ref{sec:Preliminaries}, we propose a joint estimator of the scaling parameter $\sigma$ and the Hurst parameter $H$ for the scaled Rosenblatt process (i.e. assuming $\lambda=0$ in equation \eqref{SDE}) in section~\ref{sec:scaled RP}. This estimator is strongly consistent (under the in-fill asymptotics) with a faster convergence rate than the estimators studied in \cite{AssaadTudor20} and \cite{BarTud10}. We also provide the limiting joint distribution of the estimator, being the Rosenblatt distribution. In section \ref{sec:ROU}, we address the problem of parameter identification in SDEs \eqref{SDE} that seems not to have been addressed in the literature so far. Firstly, in subsection \ref{subsec: Estimating drift}, we start with the identification of the drift parameter $\lambda$ only (with known values of the scaling parameter $\sigma$ and Hurst parameter $H$) and we construct a strongly consistent estimator under in-fill asymptotics with a fixed time horizon. The existence of such an estimator is a rather surprising result which already implies singularity of distributions of solutions corresponding to different drift parameters and discards Girsanov-type theorems for Rosenblatt process. To our best knowledge, the question about singularity/equivalence of measures generated by Rosenblatt processes with different drifts has been open so far. Secondly, in subsection \ref{subsec: all parameters}, we adapt the drift estimator (by certain decelerating procedure) and combine it with the estimator of parameters $\sigma$ and $H$ to get a jointly consistent (under in-fill asymptotics with a fixed time horizon) estimator for all three parameters. To our best knowledge, no such consistent joint estimator has been proposed yet. Finally, in section \ref{sec:Simulations} a simulation study is presented to illustrate theoretical results on convergence of estimators and their actual behaviour in finite-sample setting. 

\section{Preliminaries}\label{sec:Preliminaries}
We begin by recalling the definition of the Rosenblatt process. Let $(\Omega,\F,\P)$ be a complete probability space with a standard Wiener process $W=\{W(t)\}_{t\in [0,1]}$ defined on it. We assume that the $\sigma$-algebra $\F$ is generated by $W$. The Rosenblatt process $Z^H=\{Z^H(t)\}_{t\in [0,1]}$ with Hurst parameter $H\in (1/2,1)$ is defined here by \[Z^H(t) := I_2(L_t^H), \qquad t\in [0,1],\]
where $I_2$ denotes the Wiener-It\^o multiple integral of order $2$ with respect to the standard Wiener process $W$; see, e.g., \cite{NouPec12} or \cite{Nua06}, and where $L_t^H$ is defined by
	\[ L_t^H(x_1,x_2) := C_H^Zx_1^{-\frac{H}{2}}x_2^{-\frac{H}{2}} \left[\int_{x_1\vee x_2}^t u^H (u-x_1)^{\frac{H}{2}-1}(u-x_2)^{\frac{H}{2}-1}\d{u}\right]\bm{1}_{(0,t)^2}(x_1,x_2)\]
with the normalizing constant
	\[ C_H^Z := \frac{\sqrt{2H(2H-1)}}{2\mathrm{B}\left(1-H,\frac{H}{2}\right)};\] see, e.g., \cite{Tud08}. Here, $\mathrm{B}(\cdot,\cdot)$ denotes the Beta function. Let us also recall the 2-variation statistic for the Rosenblatt process
\begin{equation}\label{def:2variation}
V_N(Z^H):= \frac{1}{N} \sum_{i=1}^{N}\left(\frac{\left|Z^H\left(\frac{i}{N}\right) - Z^H\left(\frac{i-1}{N}\right)\right|^2}{N^{-2H}} - 1 \right),
\end{equation}
studied in \cite{TudorViens09}. In there, the authors prove that the (suitably normalized) $2$-variation $V_N(Z^H)$ converges to $Z^H(1)$ in the space $L^2(\Omega)$. It turns out, however, that the convergence also holds almost surely. More precisely, we have the following result.

\begin{Theorem}
\label{Thm:V(Z) convergence}
For any $H\in (1/2,1)$, there is the following convergence:
	\begin{equation}\label{eq:V(Z) a.s. convergence}
\frac{N^{1-H}}{4d(H)}V_N(Z^H) \quad\xrightarrow[N \to \infty]{a.s.}\quad Z^H(1).
	\end{equation}
\end{Theorem}

\begin{proof}
For $n\in\mathbb{N}_0$, denote the multiple Wiener-It\^o integral with respect to the standard Wiener process $W$ of order $n$ by $I_n$. Denote also the $n$\textsuperscript{th} Wiener chaos generated by the Wiener-It\^o integral $I_1$ by $\mathcal{H}_n$ (see, e.g., \cite{NouPec12} or \cite{Nua06} for the precise definitions). For $i=1,2,\ldots, N$, let 
	\[ A_i := L_{\frac{i}{N}}^H - L_{\frac{i-1}{N}}^H\]
and let also 
	\[ T_2 := 4N^{2H-1}\sum_{i=1}^N I_2(A_i\otimes_1A_i)\quad\mbox{and}\quad T_4 := 4N^{2H-1}\sum_{i=1}^N I_4 (A_i\otimes A_i),\]
where $\otimes_1$ denotes the $1$-contraction of indexes while $\otimes$ denotes the tensor product. If
	\[ U_2:= T_2 - \frac{4d(H)}{N^{1-H}}Z^H(1),\]
the product formula for multiple integrals from, e.g., \cite[Theorem 5.9]{CouOnd22}, provides us with the decomposition
\[
F_N:=\frac{N^{1-H}}{4d(H)}V_N(Z^H) - Z^H(1) = \frac{N^{1-H}}{4d(H)}(T_4 + U_2),
\]
see \cite[formula (3.32)]{TudorViens09}. It then follows, by \cite[Lemma 6]{TudorViens09}, that $U_2 = O_{L^2(\Omega)}(N^{-1/2})$ as $N \to \infty$, and it also follows, by the results on \cite[p. 16]{TudorViens09}, that
\[
T_4 = 
\begin{cases}
O_{L^2(\Omega)}(N^{-1/2}), &\quad H \in (\frac{1}{2},\frac{3}{4}),\\
O_{L^2(\Omega)}(\sqrt{\log N} N^{-1/2}), &\quad H = \frac{3}{4},\\
O_{L^2(\Omega)}(N^{2H-2}), &\quad H \in (\frac{3}{4},1).\\
\end{cases}
\]
Consequently, there is the following asymptotic behavior:
\begin{equation} \label{eq:2var L2 speed}
F_N=  
\begin{cases}
O_{L^2(\Omega)}(N^{\frac{1}{2}-H}), & \quad H \in (\frac{1}{2},\frac{3}{4}),\\
O_{L^2(\Omega)}(\sqrt{\log N} N^{\frac{1}{2}-H}), &\quad H = \frac{3}{4},\\
O_{L^2(\Omega)}(N^{H-1}), & \quad H \in (\frac{3}{4},1).\\
\end{cases}
\end{equation}
In all the three cases, we have that there exist $\delta>0$, $C>0$, and $N_0\in\mathbb{N}$ such that 
	\[ (\mathbb E| F_N|^2)^\frac{1}{2} \leq C N^{-\delta}\]
holds for all $N\in\mathbb{N}$, $N\geq N_0$. Let now $\kappa\in(0,\delta)$ and $\eta > \frac{1}{\delta-\kappa}$. By using Markov's inequality, we obtain the estimate
	\[ \mathbb{P}(|F_N|> N^{-\kappa}) \leq N^{\kappa\eta}\mathbb{E}|F_N|^\eta.\]
Moreover, as we clearly have that 
	\[ F_N\in \mathcal{H}_2\oplus\mathcal{H}_4,\]
the hypercontractivity result of \cite[Proposition 2.2]{CouMasOnd22} yields the existence of a constant $C_\eta>0$ such that the inequality
	\[ \mathbb E|F_N|^\eta \leq C_\eta (\mathbb E|F_N|^2)^\frac{\eta}{2}\]
holds. Thus we obtain the estimate
	\[ \mathbb{P}(|F_N|> N^{-\kappa}) \leq C_\eta C^{\eta} N^{-\eta(\delta-\kappa)},\]
from which it follows that 
	\[ \sum_{N=N_0}^\infty \mathbb{P}(|F_N|> N^{-\kappa}) <\infty.\]
The sought almost sure convergence \[F_N\xrightarrow[N \to \infty]{a.s.} 0\] is now a consequence of the Cantelli lemma.
\end{proof}

\begin{Remark}\label{rem: speed of 2var to Z(1)}
The rate of the $L^2(\Omega)$-convergence in \eqref{eq:2var L2 speed} allows to prove the convergence
\begin{equation} \label{eq:2var a.s. speed}
N^{\alpha}\left|\frac{N^{1-H}}{4d(H)}V_N(Z^H) - Z^H(1)\right| \xrightarrow[N \to \infty]{a.s.} 0
\end{equation}
for 
\[ \alpha <  
\begin{cases}
H - \frac{1}{2}, & \quad H \in (\frac{1}{2},\frac{3}{4}],\\
1-H, & \quad H \in (\frac{3}{4},1),\\
\end{cases} \]
by using the same arguments (equivalence of moments and the Cantelli lemma). This convergence rate will be useful in the sequel.
\end{Remark}
%\begin{Corollary}\label{cor:V(Z) a.s. convergence}
%For any $\gamma > 0$ we have
%\begin{equation}\label{eq:V(Z) a.s. convergence}
%N^{1-H-\gamma} V_N(Z^H) \xrightarrow[N \to \infty]{a.s.} 0,
%\end{equation}
%\end{Corollary}

%**********************************************************************************************************************

\section{Scaled Rosenblatt process}\label{sec:scaled RP}

In this section, we consider the scaled and shifted Rosenblatt process
\begin{equation}\label{eq:scaled Rosenblatt}
X(t) = X_0 + \sigma Z^H(t), \quad t \in [0,1],
\end{equation}
where $X_0\in\mathbb{R}$ is known and $H\in (1/2,1)$ as well as $\sigma>0$ are to be estimated. The following almost sure convergence is a direct consequence of \autoref{Thm:V(Z) convergence} and a variation of the $\delta$-method.

\begin{Lemma}\label{Lem: V(sZ) convergence} There is the almost sure convergence
\begin{equation}\label{Eq: V(sZ) convergence}
\frac{1}{4d(H)}N^{1-H} 
\begin{bmatrix} 
\log\left(\frac{1}{N}\sum_{i=1}^{N} \frac{|X(\frac{i}{N}) - X(\frac{i-1}{N})|^2}{\sigma^2 N^{-2H}} \right) \\
\log\left(\frac{1}{N/2}\sum_{i=1}^{N/2} \frac{|X(\frac{i}{N/2}) - X(\frac{i-1}{N/2})|^2}{\sigma^2 (N/2)^{-2H}} \right) \\
\end{bmatrix}
\quad
\xrightarrow[N \to \infty]{a.s.} 
\quad
\begin{bmatrix} 
Z^H(1)\\
2^{1-H} Z^H(1)
\end{bmatrix},
\end{equation}
where $N$ runs over even numbers.
Moreover, for any 
\[ \alpha <  
\begin{cases}
H - \frac{1}{2}, & \quad H \in (\frac{1}{2},\frac{3}{4}],\\
1-H, & \quad H \in (\frac{3}{4},1),\\
\end{cases} \]
we have the speed of convergence
\begin{equation} \label{Eq:log 2var a.s. speed}
N^{\alpha}\left| \frac{1}{4d(H)}N^{1-H} 
\begin{bmatrix} 
\log\left(\frac{1}{N}\sum_{i=1}^{N} \frac{|X(\frac{i}{N}) - X(\frac{i-1}{N})|^2}{\sigma^2 N^{-2H}} \right) \\
\log\left(\frac{1}{N/2}\sum_{i=1}^{N/2} \frac{|X(\frac{i}{N/2}) - X(\frac{i-1}{N/2})|^2}{\sigma^2 (N/2)^{-2H}} \right) \\
\end{bmatrix}
- 
\begin{bmatrix} 
Z^H(1)\\
2^{1-H} Z^H(1)
\end{bmatrix}
\right| \xrightarrow[N \to \infty]{a.s.} 0.
\end{equation}
\end{Lemma}

\begin{proof}
It is sufficient to prove the convergence of the two vector elements separately to obtain the joint convergence \eqref{Eq: V(sZ) convergence}. Moreover, it is also sufficient to only prove the convergence of the first element as the second one is its direct consequence.
Recall that \[V_N(Z^H) =\frac{1}{N} \sum_{i=1}^{N}\left(\frac{\left|Z^H\left(\frac{i}{N}\right) - Z^H\left(\frac{i-1}{N}\right)\right|^2}{N^{-2H}} - 1 \right).\] From convergence \eqref{eq:V(Z) a.s. convergence} we have that $V_N(Z^H)(\omega) = O(1/N^{1-H})$ as $N \to \infty$ for $\P$-almost all $\omega \in \Omega$. Now use the Taylor approximation 
\[ \log(x+1) = x + O(x^2), \quad \text{as} \quad x\to 0,\]
to calculate
\begin{align*}
	\frac{1}{4d(H)}N^{1-H}\log\left(\frac{1}{N}\sum_{i=1}^{N} \frac{|X(\frac{i}{N}) - X(\frac{i-1}{N})|^2}{\sigma^2 N^{-2H}} \right)
	&= \frac{1}{4d(H)}N^{1-H}\log\left(V_N(Z^H) + 1\right) \\
	&= \frac{1}{4d(H)}N^{1-H}\left(V_N(Z^H) + r_N\right).
\end{align*}
Now we have that $\frac{N^{1-H}}{4d(H)} V_N(Z^H)\xrightarrow[n\to\infty]{a.s.} Z^H(1)$ holds by \eqref{eq:V(Z) a.s. convergence}. For the second term, note that for $\mathbb{P}$-almost all $\omega\in \Omega$ we have \[r_N(\omega)= O (V_N(Z^H(\omega)))^2 =O ( O(1/N^{1-H}) )^2 = O(1/N^{2-2H}), \quad \mbox{as}\quad N \to \infty,\] and thus it follows that \[\frac{N^{1-H}}{4d(H)} r_N(\omega) =  O(1/N^{1- H}) \quad\xrightarrow[N\to\infty]{}\quad 0\] which concludes the proof of \eqref{Eq: V(sZ) convergence}. To prove convergence \eqref{Eq:log 2var a.s. speed}, it suffices to combine convergence \eqref{eq:2var a.s. speed} with the fact that for $\mathbb{P}$-almost all $\omega\in \Omega$ it holds that
\[N^{\alpha} \frac{N^{1-H}}{4d(H)} r_N(\omega) \quad\xrightarrow[N\to\infty]{}\quad 0.\]
\end{proof}

Let us now proceed to parameter estimation. The vector of logarithms on the left-hand side of \eqref{Eq: V(sZ) convergence} can be rewritten as a linear model:
\begin{align*}
\renewcommand\arraystretch{1.8}
\begin{bmatrix} 
\log\left(\frac{1}{N}\sum_{i=1}^{N} \frac{|X(\frac{i}{N}) - X(\frac{i-1}{N})|^2}{\sigma^2 N^{-2H}} \right) \\
\log\left(\frac{1}{N/2}\sum_{i=1}^{N/2} \frac{|X(\frac{i}{N/2}) - X(\frac{i-1}{N/2})|^2}{\sigma^2 (N/2)^{-2H}} \right) \\
\end{bmatrix}
\\
& 
\renewcommand\arraystretch{1.8}
\hspace{-3cm}=
\begin{bmatrix} 
(2H-1)\log N - 2 \log \sigma + \log \left(\sum_{i=1}^{N} |X(\frac{i}{N}) - X(\frac{i-1}{N})|^2\right) \\
(2H-1)\log (N/2) - 2 \log \sigma + \log \left(\sum_{i=1}^{N/2} |X(\frac{i}{N/2}) - X(\frac{i-1}{N/2})|^2\right) \\
\end{bmatrix}
\\
&
\renewcommand\arraystretch{1.8}
\hspace{-3cm}=
\underbrace{\begin{bmatrix} 
\log N 			& 1\\
\log (N/2)  & 1\\
\end{bmatrix}}_{=:\mathbb{X}}
\cdot
\underbrace{\begin{bmatrix} 
 2H - 1 \\
-2 \log \sigma\\
\end{bmatrix}}_{=:\beta}
+
\underbrace{\begin{bmatrix} 
 \log \left(\sum_{i=1}^{N} |X(\frac{i}{N}) - X(\frac{i-1}{N})|^2\right) \\
 \log \left(\sum_{i=1}^{N/2} |X(\frac{i}{N/2}) - X(\frac{i-1}{N/2})|^2\right)\\
\end{bmatrix}}_{=:\mathbb{Y}}.
 \numberthis\label{Eq: V(sZ) linear model} 
\end{align*}
As such, in order to separate vector $\beta$, that contains the parameters of interest, we can use the least-squares formalism. By \eqref{Eq: V(sZ) convergence}, we have 
	\[ \frac{N^{1-H}}{4d(H)} (\mathbb{X} \cdot \beta + \mathbb{Y}) \quad\xrightarrow[\substack{N\in 2\mathbb{N}\\ N\to\infty}]{a.s.}\quad \begin{bmatrix} Z^H(1) \\ 2^{1-H} Z^H(1) \end{bmatrix}.\]
We also have that 
\begin{equation}\label{eq:LSE factor}
W_N :=
\begin{bmatrix} 
  1 			& 0\\
  0       & \frac{1}{\log N} \\
\end{bmatrix}
\cdot
(\mathbb{X}^T \cdot \mathbb{X})^{-1} \cdot \mathbb{X}^T = \frac{1}{\log 2} 
\begin{bmatrix} 
  1 											      & -1\\
  -1 + \frac{\log 2}{\log N}    & 1 \\
\end{bmatrix}
\quad \xrightarrow[N \to \infty]{}\quad  \frac{1}{\log 2} 
\begin{bmatrix} 
  1 		& -1\\
  -1    & 1 \\
\end{bmatrix}.
\end{equation}
It thus follows, that 
	\[ W_N \cdot \frac{N^{1-H}}{4d(H)} (\mathbb{X}\cdot\beta + \mathbb{Y}) \quad \xrightarrow[\substack{N\in 2\mathbb{N}\\ N\to\infty}]{a.s.}\quad  \begin{bmatrix} -1 \\ 1 \end{bmatrix} \cdot \frac{2^{1-H} - 1 }{\log 2} Z^H(1).\]
On the other hand, if we set 
	\[
\hat{\beta}_N :=  -(\mathbb{X}^T \cdot \mathbb{X})^{-1} \cdot \mathbb{X}^T\mathbb{Y},
\]
we also have
\begin{equation*}
W_N \cdot \frac{N^{1-H}}{4d(H)} (\mathbb{X} \cdot \beta + \mathbb{Y})   = 
\frac{N^{1-H}}{4d(H)} \begin{bmatrix} 
  1 			& 0\\
  0       & \frac{1}{\log N} \\
\end{bmatrix}
\cdot
\left( \beta -\hat\beta_N\right)
\end{equation*}
so that 
	\begin{equation}
	\label{eq:beta convergence}
	 \frac{N^{1-H}}{4d(H)} \begin{bmatrix} 
  1 			& 0\\
  0       & \frac{1}{\log N} \\
\end{bmatrix}
\cdot
\left( \beta -\hat\beta_N\right) \quad\xrightarrow[\substack{N\in 2\mathbb{N}\\ N\to\infty}]{a.s.}\quad \begin{bmatrix} -1 \\ 1 \end{bmatrix} \cdot \frac{2^{1-H} - 1 }{\log 2} Z^H(1).
\end{equation}
This enables us to express estimators for $H$ and $\sigma$ explicitly and, moreover,  we obtain the speed of the almost sure convergence together with the limiting random variable. These considerations result in the definition and strong consistency of the estimators of parameters $H$ and $\sigma$. We formulate this as a theorem.

\begin{Theorem}\label{thm: estimators convergence for sRP} 
Consider the following two estimators:
\begin{equation}\label{def: estimators H sigma}
\begin{aligned}
\hat{H}_N & := - \frac{1}{2 \log 2}\Bigg[ \log \left(\sum_{i=1}^{N} \left|X\left(\frac{i}{N}\right) - X\left(\frac{i-1}{N}\right)\right|^2\right) \\
& \hspace{3cm} - \log \Bigg(\sum_{i=1}^{N/2} \left|X\left(\frac{i}{N/2}\right) - X\left(\frac{i-1}{N/2}\right)\right|^2\Bigg) \Bigg] + \frac{1}{2},\\
\hat{\sigma}_N & := \exp 
	\Bigg\{\frac{1}{2 \log 2}  
		\Bigg[\log\left(\frac{2}{N}\right) \log \left(\sum_{i=1}^{N} \left|X\left(\frac{i}{N}\right) - X\left(\frac{i-1}{N}\right)\right|^2\right) \\
		& \hspace{3cm} + \log(N)  \log \Bigg(\sum_{i=1}^{N/2} \left|X\left(\frac{i}{N/2}\right) - X\left(\frac{i-1}{N/2}\right)\right|^2\Bigg) 
		\Bigg] 
	\Bigg\} 
\end{aligned}
\end{equation}
These estimators are strongly consistent with the following speed of convergence (and the limit)
\begin{equation}\label{eq: convergence H sigma}
\begin{bmatrix} 
 N^{1-H} (\hat{H}_N - H) \\
 \frac{N^{1-H}}{\log N} (\hat{\sigma}_N - \sigma) \\
\end{bmatrix} 
\quad \xrightarrow[\substack{N\in 2\mathbb{N}\\ N\to\infty}]{a.s.} \quad
\begin{bmatrix} 
1\\
\sigma\\
\end{bmatrix}
\frac{2 d(H) (2^{1-H}-1)}{\log 2} Z^H(1).
\end{equation}
\end{Theorem}

\begin{proof}
The claim follows directly by the convergence in \eqref{eq:beta convergence} and by the application of the $\delta$-method similarly as in the proof of \autoref{Lem: V(sZ) convergence}.
\end{proof}

\section{SDEs driven by additive Rosenblatt process}\label{sec:ROU}

Let $F:\R\to\R$ be a Borel measurable function such that
\begin{enumerate}[label=(F\arabic*)]
	\item\label{F1} $F$ is locally Lipschitz, i.e.\ for every $N\in\mathbb{N}$ there exists a constant $K_N\in (0,\infty)$ such that for every $x,y\in\mathbb{R}$, $|x|+|y|\leq N$, it holds that \[|F(x)-F(y)|\leq K_N |x-y|,\] 
	\item\label{F2} and $F$ satisfies a Lyapunov condition:\ There exists a function $V\in C^1(\R)$ which satisfies
		\begin{equation}
		\label{F2a} \lim_{R\to\infty} \inf_{|x|>R} V(x)=\infty
		\end{equation}
 and for which there is a constant $K\in (0,\infty)$ and a continuous function $h:[0,\infty)\to [0,\infty)$ such that the inequality
		\begin{equation}
		\label{F2b}
		 V'(x)F(x+z) \leq K V(x) + h(|z|)
		\end{equation}
	holds for every $x,z\in\R$. 
\end{enumerate}
Let us remark that condition \ref{F2} relaxes the usual linear growth condition. Indeed, by setting $V(x)=x^2$, we obtain that $F$ is required to satisfy $xF(x) \leq k_1 + k_2x^2$ for some constants $k_1,k_2\in (0,\infty)$. This, in particular, allows for $F$ to be a polynomial of an odd order with a negative leading coefficient. Let now $Z=(Z_t,t\in [0,1])$ be a stochastic process defined on some probability space $(\Omega,\mathcal{F},\mathbb{P})$ that has $\nu$-H\"older continuous sample paths for $\nu\in (0,1)$ and let $\omega\in\Omega$ be such that this holds for the path $Z(\omega)$. Then, by condition \ref{F1}, the map $(t,x)\mapsto F(x+Z_t(\omega))$ satisfies Carath\'eorody's conditions so that the equation 
	\[ y_\omega^\prime(t) = F(y_\omega(t) + Z_t(\omega))\]
subject to $y_{\omega}(0) = x_0$ admits a solution $y_\omega$ defined on its maximal interval of existence $[0,\tau)$.
If one defines $w(t,x):=\mathrm{e}^{-Kt}V(x)$ for $(t,x)\in [0,1]\times\R$, we obtain, by condition \eqref{F2b}, the estimate
	\[ \frac{\d}{\d t} w(t,y_\omega(t)) \leq \mathrm{e}^{-Kt}h(|Z_t(\omega)|),\quad t\in (0,\tau),\]
and, by integrating this inequality, we also obtain 
	\[ w(t,y_\omega(t)) - \mathrm{e}^{-Kt_0}V(x_0) \leq \sup_{s\in [0,1]} h(|Z_s(\omega)|), \qquad t\in [0,\tau], \]
which means that 
	\[ V(y_\omega(t)) \leq \mathrm{e}^{K}\left(V(x_0) + \sup_{s\in [0,1]} h(|Z_s(\omega)|)\right), \qquad t\in [0,\tau],\]
where the right-hand side is finite by continuity of $h$ and $Z(\omega)$. But $V$ satisfies \eqref{F2a} and it therefore follows that $y_\omega$ must be bounded on $[0,\tau]$. In particular, it does not explode at $t=\tau$, and hence, we must have $\tau=1$.
Define now $X(\omega):=(X_t(\omega), t\in [0,1])$ by
	\[ X_t(\omega):= y_\omega(t) + Z_t(\omega), \qquad t\in [0,1].\]
Then $X(\omega)$ satisfies the differential equation
	\[ X_t(\omega) = x_0 + \int_0^t F(X_s(\omega))\d{s} + Z_t(\omega), \qquad t\in [0,1]. \]
Moreover, we have that 
	\[\sup_{t\in [0,1]}|X_t(\omega)| \leq \sup_{t\in [0,1]}|y_\omega(t)| + \sup_{t\in [0,1]} |Z_t(\omega)| = :C(\omega), \]
where $C(\omega)$ is finite by boundedness of $y_\omega$ proved above and by continuity of $Z(\omega)$, and also 
	\begin{align*}
	 |X_t(\omega)-X_s(\omega)| & \leq |y_\omega(t)-y_\omega(s)| + |Z_t(\omega)-Z_s(\omega)|\\
	 &  \leq \int_s^t |F(X_r(\omega))|\d{r} + \|Z(\omega)\|_{C^\nu([0,1])} \,|t-s|^\nu \\
	& \leq K_{\lceil C(\omega)\rceil}(C(\omega) + |F(0)|)\,|t-s| + \|Z(\omega)\|_{C^\nu([0,1])} \,|t-s|^\nu
	 \end{align*}
for $s,t\in [0,1]$ so that $X(\omega)$ is $\nu$-H\"older continuous on $[0,1]$.
%Let $F:\mathbb{R}\to\mathbb{R}$ be a Borel measurable function that satisfies the following two conditions:
%	\begin{itemize}
%		\item $F$ has at most linear growth, i.e.\ there exists a constant $K_F\in (0,\infty)$ such that for every $x\in\mathbb{R}$ it holds that \[|F(x)| \leq K_F(1+|x|);\] 
%		\item and $F$ is locally Lipschitz, i.e.\ for every $N\in\mathbb{N}$ there exists a constant $K_N\in (0,\infty)$ such that for every $x,y\in\mathbb{R}$, $|x|+|y|\leq N$, it holds that \[|F(x)-F(y)|\leq K_N |x-y|.\] 
%	\end{itemize} 
%Let also $Z$ be a continuous stochastic process, defined on some probability space $(\Omega,\mathcal{F},\mathbb{P})$, and $X_0\in\mathbb{R}$. It follows, by a standard Picard iteration scheme, that there is a unique (in the sense of indistinguishability) continuous stochastic process $X$ that satisfies the random differential equation 
%	\[ X(t) = X_0 + \int_0^t F(X_s) \d{s} + Z_t, \quad t\in [0,1],\, a.s.;\]
%cf., e.g., Proposition 3.5 in \cite{Snup09}. Moreover, it can be shown by the same argument as in Proposition 3.6 in \cite{Snup09} that if there exists $\nu\in (0,1)$ such that the stochastic process $Z$ has a version with sample paths in $C^\gamma([0,1])$ for every $\gamma\in (0,\nu)$ almost surely, then the solution $X$ also has a version with sample paths in the H\"older space $C^\gamma([0,1])$ for every $\gamma\in (0,\nu)$ almost surely and the estimate 
%	\[ \|X\|_{C^\gamma([0,1])} \leq C (1+ \|Z\|_{C^\gamma([0,1])})\]
%holds almost surely with a finite positive constant $C\equiv C(X_0,K_F,\gamma)$.
By applying these results to the scaled Rosenblatt process $Z=\sigma Z^H$ (with $\sigma>0$ and $H\in (1/2,1)$) and $F=\lambda f$ (with $\lambda\in\mathbb{R}$ and $f$ satisfying the local Lipschitz condition \ref{F1} and the Lyapunov condition \ref{F2}), we obtain a unique process $X$ that satisfies 
	\begin{equation}
	\label{eq:non-lin_SDE}
	 X(t) = X_0 + \lambda \int_0^t f(X(s))\d{s} + \sigma Z^H(t), \qquad t\in [0,1],\, a.s.,
	\end{equation}
with sample paths in $C^{\gamma}([0,1])$ for every $\gamma\in (0,H)$ almost surely.
%In this section the solution $X$ to the linear SDE \eqref{SDE} with non-trivial drift is studied. Note that 
%\begin{equation}\label{eq:ROU reseni SDE}
%X(t) = X_0 + \lambda \int_0^t X(s) \d s + \sigma Z^H(t) 
%\end{equation}
%holds for $t\in [0,1]$ and that 
We can decompose the corresponding 2-variation as follows:
\begin{align*}
\frac{N^{1-H}}{4 d(H)} V_N(X) & = \frac{N^{1-H}}{4 d(H)}\frac{1}{N} \sum_{i=1}^{N}\left(\frac{\left|X\left(\frac{i}{N}\right) - X\left(\frac{i-1}{N}\right)\right|^2}{\sigma^2 N^{-2H}} - 1 \right) \\
& = \frac{N^{1-H}}{4 d(H)}\frac{1}{N} \sum_{i=1}^{N}\frac{\lambda^2}{\sigma^2 N^{-2H}} \left(\int_{\frac{i-1}{N}}^{\frac{i}{N}} f(X(s)) \d s\right)^2  \\
&\hspace{1cm} +\frac{N^{1-H}}{2 d(H)} \frac{1}{N} \sum_{i=1}^{N}\frac{\lambda}{\sigma N^{-2H}}  \left(\int_{\frac{i-1}{N}}^{\frac{i}{N}} f(X(s)) \d s\right) \left(Z^H\left(\frac{i}{N}\right) - Z^H\left(\frac{i-1}{N}\right)\right) \\
&\hspace{1cm} +\frac{N^{1-H}}{4 d(H)} \frac{1}{N} \sum_{i=1}^{N}\left(\frac{1}{ N^{-2H}}\left(Z^H\left(\frac{i}{N}\right) - Z^H\left(\frac{i-1}{N}\right)\right)^2 - 1 \right) \\
&	 =: A_N + B_N + C_N.	\numberthis\label{eq:2var X}
\end{align*}
Let us now calculate the limits of each of the terms separately. We start with term $B_N$ and show that it converges almost surely to zero. Fix an $\omega \in \Omega$ and write
\begin{align*}\label{eq:lim B_N}
B_N(\omega)  & = \frac{\lambda}{2 d(H) \sigma} N^H  \sum_{i=1}^{N} \left(\int_{\frac{i-1}{N}}^{\frac{i}{N}} f(X(s)(\omega)) \d s\right) \left(Z^H\left(\frac{i}{N}\right)(\omega) - Z^H\left(\frac{i-1}{N}\right)(\omega)\right) \\
& = \frac{\lambda}{2 d(H) \sigma} N^{H-1}  \sum_{i=1}^{N} f(X(s_i)(\omega)) \left(Z^H\left(\frac{i}{N}\right)(\omega) - Z^H\left(\frac{i-1}{N}\right)(\omega)\right),
\end{align*}
for some $(i-1)/N \leq s_i \leq i/N$ that can depend on $\omega$. We also have that for every $\gamma<H$, there exists $M(\omega)\in\mathbb{N}$ such that  
	\[ |X(t)(\omega)| + |X(s)(\omega) | \leq 2\|X(\cdot)(\omega)\|_{C^\gamma ([0,1])} \leq M(\omega)\]
holds for every $s,t\in [0,1]$, $s<t$, and therefore
	\begin{equation}
	\label{eq:f(X)_Holder}
	 |f(X(t)(\omega))-f(X(s)(\omega))| \leq K_{M(\omega)} |X(t)(\omega)-X(s)(\omega)| \leq K_{M(\omega)} \| X(\cdot)(\omega)\|_{C^{\gamma}([0,1])} |t-s|^\gamma
	 \end{equation}
holds for any $s,t\in [0,1]$, $s<t$, by appealing to the local Lipschitz condition for function $f$ and $\gamma$-H\"older continuity of $X(\cdot)(\omega)$. It thus follows that the path $f(X(\cdot)(\omega))$ is $\gamma$-H\"older continuous. As we also have that $Z^H(\cdot)(\omega)$ is $\gamma$-H\"older continuous and $H>1/2$, we can employ, path by path, the celebrated Young's result on existence of the Stieltjes ingeral in the Riemann sense (cf. \cite{Young36}) to get the convergence:
\begin{equation}\label{eq:pathwise integral}
\sum_{i=1}^{N} f(X(s_i)(\omega)) \left(Z^H\left(\frac{i}{N}\right)(\omega) - Z^H\left(\frac{i-1}{N}\right)(\omega)\right) \xrightarrow[N \to \infty]{}  \int_0^1 f(X(s)(\omega)) \circ \d Z^H(s)(\omega) < \infty.
\end{equation}
Thus, we obtain the almost sure convergence:
\[B_N \xrightarrow[N \to \infty]{a.s.} 0.\]
The term $A_N$ can be analyzed in a similar manner with the difference that the increments of $Z^H$ are replaced by the increments of a more regular (continuously differentiable) indefinite Riemann integral $\int_0^t f(X(s)) \d s$. This leads to the almost sure convergence:
\[A_N \xrightarrow[N \to \infty]{a.s.} 0.\]
Finally, the almost sure convergence of $C_N$ follows directly from \autoref{Thm:V(Z) convergence}:
\[
C_N =  \frac{N^{1-H}}{4d(H)}V_N(Z^H) \xrightarrow[N \to \infty]{a.s.} Z^H(1).
\]
Hence, we obtain the following convergence result for 2-variations of the solution to the SDE \eqref{eq:non-lin_SDE}:

\begin{Theorem}\label{Thm:V(X) convergence}
The normalized 2-variation of the solution $X$ to equation \eqref{eq:non-lin_SDE} converges almost surely to the actual realization of $Z^H$ at time $t=1$:
\begin{equation}\label{eq:V(X) convergence}
\frac{N^{1-H}}{4 d(H)}\frac{1}{N} \sum_{i=1}^{N}\left(\frac{\left|X\left(\frac{i}{N}\right) - X\left(\frac{i-1}{N}\right)\right|^2}{\sigma^2 N^{-2H}} - 1 \right) \quad \xrightarrow[N \to \infty]{a.s.}\quad Z^H(1).
\end{equation}
%Moreover, for any $\gamma > 0$ we have
%\begin{equation}\label{eq:V(X) a.s. convergence}
%\frac{N^{1-H-\gamma}}{4 d(H)}\frac{1}{N} \sum_{i=1}^{N}\left(\frac{|X(i/N) - X((i-1)/N)|^2}{\sigma^2 N^{-2H}} - 1 \right)  \xrightarrow[N \to \infty]{a.s. } 0.
%\end{equation}
\end{Theorem}

\begin{Remark}\label{diff estim for ROU}
We now have all that is needed to repeat the ideas and calculations from \autoref{sec:scaled RP}. This results in the fact that all the properties of the estimators $\hat{H}_N$ and $\hat{\sigma}_N$ given in \autoref{thm: estimators convergence for sRP} remain valid also for the solution $X$ to equation \eqref{eq:non-lin_SDE}.
\end{Remark}

\begin{Remark}
In \cite{AssaadTudor20}, the convergence of the 2-variation of an Hermite Ornstein--Uhlenbek process of type \eqref{eq:V(X) convergence} is studied as well but the convergence is understood in the $L^2(\Omega)$ sense. Compared to the (relatively simple) pathwise argument for the almost sure convergence of term $B_N$ presented above, the calculations in \cite{AssaadTudor20}, that are needed to prove the $L^2(\Omega)$ convergence, are rather complicated. In addition, the pathwise argument given here directly generalizes to solutions of non-linear equations provided that their solutions are sufficiently smooth.   
\end{Remark}

\subsection{Estimating drift}\label{subsec: Estimating drift}
For the purposes of this subsection, we assume that parameters $H$ and $\sigma$ are known. We show that in this situation, it is possible to find a consistent estimate of the drift parameter of the solution to SDE \eqref{eq:non-lin_SDE} from a discretely observed single trajectory with a fixed time horizon and a decreasing time mesh. We note that this is in sharp contrast with the case of SDEs driven by a (fractional) Brownian motion where such estimation is impossible because the probability laws of the solutions corresponding to different drift parameters on the space of continuous trajectories are equivalent (by the Girsanov theorem for fractional Brownian motions established in \cite{DecrUstu99}).

If parameters $H$ and $\sigma$ are known, the 2-variation which approximates values of the driving process $Z^H$ well can be evaluated. In particular, notice that \autoref{Thm:V(Z) convergence} and self-similarity and stationarity of the increments of the Rosenblatt process imply the following two convergences:

\begin{equation}\label{eq: V(Z) to Z1/2}
\frac{N^{1-H}}{4d(H)} \frac{1}{N} \sum_{i=1}^{N/2}\left(\frac{\left|Z^H\left(\frac{i}{N}\right) - Z^H\left(\frac{i-1}{N}\right)\right|^2}{N^{-2H}} - 1 \right) \quad\xrightarrow[\substack{N\in 2\mathbb{N}\\ N\to\infty}]{a.s.}\quad Z^H(1/2),
\end{equation} 

\begin{equation}\label{eq: V(Z) to Z1 - Z1/2}
\frac{N^{1-H}}{4d(H)} \frac{1}{N} \sum_{i=N/2+1}^{N}\left(\frac{\left|Z^H\left(\frac{i}{N}\right) - Z^H\left(\frac{i-1}{N}\right)\right|^2}{N^{-2H}} - 1 \right) \quad\xrightarrow[\substack{N\in 2\mathbb{N}\\ N\to\infty}]{a.s.}\quad  Z^H(1) - Z^H(1/2).
\end{equation} 
Indeed, $H$ self-similarity of the Rosenblatt process implies equality of laws
\[
\text{Law}\left(\{Z^H(t)\}_{t\geq 0}\right) = \text{Law}\left(\{2^{-H}Z^H(2t)\}_{t\geq 0}\right)
\]
on the trace of the product of Borel $\sigma$-algebras $[\mathcal{B}(\mathbb{R})]^{[0,\infty)}$ on continuous functions $C([0,\infty))$. Hence,
\begin{equation}
\begin{aligned}
&\mathbb{P}\left[  \lim_{N \to \infty}  \frac{N^{1-H}}{4d(H)} \frac{1}{N} \sum_{i=1}^{N/2} \left( \frac{\left|Z^H\left(\frac{i}{N}\right) - Z^H\left(\frac{i-1}{N}\right)\right|^2}{N^{-2H}} - 1\right) - Z^H\left(\frac{1}{2}\right) = 0 \right] \\
& = \mathbb{P}\left[ 2^{-H} \Bigg(\lim_{N \to \infty}  \frac{(N/2)^{1-H}}{4d(H)} \frac{1}{N/2} \sum_{i=1}^{N/2} \left( \frac{\left|Z^H\left(\frac{i}{N/2}\right) - Z^H\left(\frac{i-1}{N/2}\right)\right|^2}{(N/2)^{-2H}} - 1\right) - Z^H\left(1\right) \Bigg) = 0 \right] = 1,
\end{aligned}
\end{equation}
where the limits are taken over even $N$ and where the last equality follows by \autoref{Thm:V(Z) convergence}. Convergence \eqref{eq: V(Z) to Z1 - Z1/2} follows by same argument combined with stationarity of increments of $Z^H$.\\

As a consequence, convergence results for the solution process $X$ corresponding to \autoref{Thm:V(X) convergence} can be obtained by the same arguments:

\begin{equation}\label{eq:V(X1/2)  convergence}
W^{[0,1/2]}_{1/N}(X):=\frac{N^{1-H}}{4 d(H)}\frac{1}{N} \sum_{i=1}^{N/2}\left(\frac{\left|X\left(\frac{i}{N}\right) - X\left(\frac{i-1}{N}\right)\right|^2}{\sigma^2 N^{-2H}} - 1 \right)  \quad\xrightarrow[\substack{N\in 2\mathbb{N}\\ N\to\infty}]{a.s. }\quad Z^H(1/2),
\end{equation}
and 
\begin{equation}\label{eq:V(X1 - X1/2)  convergence}
W^{[1/2,1]}_{1/N}(X):=\frac{N^{1-H}}{4 d(H)}\frac{1}{N} \sum_{i=N/2+1}^{N}\left(\frac{\left|X\left(\frac{i}{N}\right) - X\left(\frac{i-1}{N}\right)\right|^2}{\sigma^2 N^{-2H}} - 1 \right)  \xrightarrow[\substack{N\in 2\mathbb{N}\\ N\to\infty}]{a.s.} Z^H(1) - Z^H(1/2).
\end{equation}
The above two approximations of the Rosenblatt path on intervals $[0,1/2]$ and $[1/2,1]$ can now be used for the following linear model:
\begin{equation}\label{Eq: linear model lambda} 
\underbrace{\begin{bmatrix} 
X(1/2) - X(0) - \sigma W^{[0,1/2]}_{1/N}(X)\\
X(1) - X(1/2) - \sigma W^{[1/2,1]}_{1/N}(X)
\end{bmatrix}}_{=:\mathbb{U}_N}
= \lambda
\underbrace{\begin{bmatrix} 
\int_{0}^{1/2}f(X(s)) \d s\\
\int_{1/2}^{1}f(X(s))\d s\\
\end{bmatrix}}_{=:\mathbb{Z}} 
+
\sigma\underbrace{\begin{bmatrix} 
Z^H(1/2) - W^{[0,1/2]}_{1/N}(X)  \\
(Z^H(1) - Z^H(1/2)) - W^{[1/2,1]}_{1/N}(X)\\
\end{bmatrix}}_{=:\varepsilon_N}.
\end{equation}
Now, it follows by \eqref{eq:V(X1/2)  convergence} and \eqref{eq:V(X1 - X1/2)  convergence} that the error term converges to zero:
\[
\mathbf{\varepsilon}_N 
\quad \xrightarrow[\substack{ N\in 2\mathbb{N}\\ N \to \infty}]{a.s.} \quad
\begin{bmatrix} 
0 \\
0 \\
\end{bmatrix}.
\]
Vector $\mathbb{Z}$ can be approximated by 
\[
\renewcommand{\arraystretch}{1.8}
\mathbb{Z}_N:= \begin{bmatrix} 
 \sum_{i=1}^{N/2}f\left(X\left(\frac{i}{N}\right)\right) \frac{1}{N} \\
 \sum_{i=N/2+1}^{N} f\left(X\left(\frac{i}{N}\right)\right) \frac{1}{N} \\
\end{bmatrix}
\]
because of the convergence $\mathbb{Z}_N\xrightarrow[]{a.s.} \mathbb{Z}$ as $N\in2\mathbb{N}$, $N\to\infty$, that follows by the continuity of sample paths of $X$ and the continuity of $f$. Both $\mathbb{U}_N$ and $\mathbb{Z}_N$ can be calculated from the available data (recall that $H$ and $\sigma$ are assumed to be known in this subsection) which leads to the following definition of the estimator of drift $\lambda$, denoted by $\hat{\lambda}_{1,N}$, by using the least-squares method:
\begin{equation}\label{def: hat(lambda)1}
\hat{\lambda}_{1,N} := (\mathbb{Z}_N^T \cdot \mathbb{Z}_N)^{-1} \cdot \mathbb{Z}_N^T \cdot \mathbb{U}_N.
\end{equation}

\begin{Theorem}\label{thm: hat(lambda)1 convergence}
The estimator $\hat{\lambda}_{1,N}$ defined by formula \eqref{def: hat(lambda)1} is strongly consistent, i.e.
\begin{equation}\label{eq: hat(lambda)1 consistent}
\hat{\lambda}_{1,N} \xrightarrow[\substack{N\in 2\mathbb{N}\\ N \to \infty}]{a.s.} \lambda,
\end{equation}
with the convergence rate at least
\begin{equation}\label{eq: hat(lambda)1 convergence}
N^{\alpha}|\hat{\lambda}_{1,N} - \lambda| \xrightarrow[\substack{N\in 2\mathbb{N}\\ N \to \infty}]{a.s.} 0
\end{equation}
for
	\[ \alpha <  
\begin{cases}
H - \frac{1}{2}, &\quad H \in (\frac{1}{2},\frac{3}{4}],\\
1-H, &\quad H \in (\frac{3}{4},1).\\
\end{cases}\]
\end{Theorem}

\begin{proof}
Recall that the convergence in \eqref{eq:2var a.s. speed} guarantees $N^{\alpha} \mathbf{\varepsilon}_N \xrightarrow[N \to \infty]{a.s.} 0$. Recall also, that it follows from the proof of \autoref{Thm:V(X) convergence} (see, in particular, inequality \eqref{eq:f(X)_Holder}) that for almost every $\omega\in\Omega$, the path $f(X(\cdot)(\omega))$ is $\gamma$-H\"older continuous for every $\gamma<H$. Hence, we have that
\begin{align*}
\left|\int_{0}^{\frac{1}{2}}f(X(s)(\omega)) \d s - \sum_{i=1}^{\frac{N}{2}}f\left(X\left(\frac{i}{N}\right)(\omega)\right)\frac{1}{N}  \right| & \leq  \sum_{i=1}^{\frac{N}{2}} \int_{\frac{i-1}{N}}^{\frac{i}{N}} \left|f(X(s)(\omega) - f\left(X\left(\frac{i}{N}\right)(\omega)\right) \right| \d s \\
& \leq \sum_{i=1}^{\frac{N}{2}} \int_{\frac{i-1}{N}}^{\frac{i}{N}} C(\omega) \left|s - \frac{i}{N} \right|^{\gamma} \d s  \\
& \leq C(\omega) N^{-\gamma}, 
\end{align*}
holds for every $\gamma<H$ with an (almost surely) finite and positive random constant $C$ that can change from line to line. It follows directly that
\begin{equation}\label{eq:Z_N convergence rate}
\left| \mathbb{Z}_N(\omega) - \mathbb{Z}(\omega) \right| \leq C(\omega) {N^{-\gamma}}
\end{equation}
holds for every $\gamma<H$ and, consequently, we have the convergence
\[
(\mathbb{Z}_N^T \cdot \mathbb{Z}_N)^{-1} \cdot \mathbb{Z}_N^T \xrightarrow[\substack{N\in 2\mathbb{N}\\ N \to \infty}]{a.s.} (\mathbb{Z}^T \cdot \mathbb{Z})^{-1} \cdot \mathbb{Z}^T.
\]
Now we can conclude the argument by computing
\begin{align*}
N^{\alpha}|\hat{\lambda}_{1,N} - \lambda| & = N^{\alpha} \left| (\mathbb{Z}_N^T \cdot \mathbb{Z}_N)^{-1} \cdot \mathbb{Z}_N^T \cdot \mathbb{U}_N  - (\mathbb{Z}_N^T \cdot \mathbb{Z}_N)^{-1} \cdot \mathbb{Z}_N^T \cdot \mathbb{Z}_N \lambda \right| \\
& \leq N^{\alpha} \left| (\mathbb{Z}_N^T \cdot \mathbb{Z}_N)^{-1} \cdot \mathbb{Z}_N^T \cdot (\mathbb{Z}-\mathbb{Z}_N) \lambda \right| +  N^{\alpha} \left| (\mathbb{Z}_N^T \cdot \mathbb{Z}_N)^{-1} \cdot \mathbb{Z}_N^T \cdot \mathbf{\varepsilon}_N \right|.
\end{align*}
and noting that the two summands above converge to zero by previous considerations (we can choose $\gamma$ in \eqref{eq:Z_N convergence rate} so that $\alpha < \gamma$). The strong consistency follows immediately.
\end{proof}

\begin{Remark} Let us note that it is shown in \cite{ChronopoulouEtAl11} that the (renormalized) quadratic variations of Hermite processes of order higher than 2 do not converge to the value of the given Hermite process at time 1 (they again converge to a Rosenblatt random variable). Thus, the estimation procedure for the drift parameter that is described above is directly applicable only for the Hermite process of order 2 (that is the Roseblatt process). The possibility to consistently estimate drift parameter for Hermite processes of higher order under fixed time-horizon regime remains (to our best knowledge) open.
\end{Remark}

\begin{Remark}\label{rem:singularity}
\textbf{(Singularity of measures)} As already noted, due to the Girsanov theorem from \cite{DecrUstu99}, a consistent estimation of the drift in SDEs driven by FBMs in the high-frequency discrete setting with a fixed time horizon is impossible (unless the time horizon tends to infinity). This suggests that a Girsanov-type theorem can not hold true in general for the (drifted) Rosenblatt process. Indeed, the estimator $\hat{\lambda}_{1,N}$ can be understood as a measurable mapping on the space $C([0,1])$ (where the trajectories of process $X$ live) and therefore, if we interpret $\hat\lambda_{1,N}$ in this manner, the mapping
\[
G(x):= \limsup_{\substack{N\in 2\mathbb{N}\\ N \to \infty}} \hat{\lambda}_{1,N}(x), \qquad x\in C([0,1]).
\]
is also measurable. 
Now consider again a Rosenblatt process $\{Z^H(t)\}_{t \in [0,1]}$ and a drifted Rosenblatt process $\{X(t) = \lambda t + Z^H(t)\}_{t\in [0,1]}$ with $\lambda \neq 0$. By the strong consistency of the estimator $\hat{\lambda}_{1,N}$, we have that both equations
\begin{equation*}
\begin{aligned}
&G(Z^H(\cdot)(\omega)) = 0,\\
&G(X(\cdot)(\omega)) = \lambda, 
\end{aligned}
\end{equation*}
hold for $\mathbb P$-almost all $\omega\in\Omega$. This, however, means that distributions of processes $Z^H$ and $X$ on the space $C([0,1])$ are mutually singular. This is in contrast with the case of (fractional) Brownian motion and the drifted (fractional) Brownian motion in which the Girsanov theorem implies that these processes have equivalent distributions on the space of trajectories (cf. \cite{DecrUstu99}).
\end{Remark}

\subsection{Estimating all parameters}\label{subsec: all parameters}
In this section, we explore the possibility of estimation of all the three parameters $\lambda, \sigma$ and $H$ in the SDE \eqref{SDE} based on high-frequency data; i.e.\ from a discretely observed single trajectory with a fixed and finite time horizon and a decreasing time mesh. 

First, by appealing to \autoref{diff estim for ROU}, we can estimate the diffusion parameters $\sigma$ and $H$ by computing the estimators $\hat{\sigma}_N$ and $\hat{H}_N$ that are given in \eqref{def: estimators H sigma}, respectively. These estimators do not require the knowledge of the drift $\lambda$. 
Second, we can try to estimate the drift parameter $\lambda$ by computing the estimator $\hat{\lambda}_{1,N}$ defined in \eqref{def: hat(lambda)1} and replace the unknown parameters $\sigma$ and $H$ in $\mathbb{U}_N$ with their estimates $\hat{\sigma}_N$ and $\hat{H}_N$, respectively. In what follows, we examine convergence of this plug-in drift estimator.

The key ingredient for the convergence of the drift estimator is the convergence of the rescaled 2-variation towards the corresponding realization of the driving process that is stated in \eqref{eq:V(X1/2)  convergence} and \eqref{eq:V(X1 - X1/2)  convergence}. For simplicity of the exposition, we study the 2-variation on the whole interval $[0,1]$ instead of on the two subintervals $[0,1/2]$ and $[1/2,1]$. Their properties are analogous. Moreover, we switch to the logarithmic version of the 2-variation, introduced in \autoref{Lem: V(sZ) convergence}, as it is more convenient for our purposes.

We now show that the convergence in \eqref{Eq: V(sZ) convergence} does not hold if we replace $H$ and $\sigma$ with their estimates. Indeed, let us set
	\begin{align*}
		m_N & := \frac{d(H)}{d(\hat{H}_N)} \e^{(H-\hat{H}_N) \log N},\\
		L_1(N) & := \frac{ N^{1-H}}{4d(H)} \left[ 2(\hat{H}_N - H) \log N  -2(\log \hat{\sigma}_N - \log \sigma) \right], \\
		L_2(N) & := \frac{ N^{1-H}}{4d(H)} \left[(2H - 1) \log N  -2 \log \sigma + \log\left(\sum_{i=1}^{N} \left|X\left(\frac{i}{N}\right) - X\left(\frac{i-1}{N}\right)\right|^2\right) \right]
	\end{align*}
and, with this notation, write
	\begin{equation}
	\label{2var non-convergence}
	 \frac{N^{1-\hat{H}_N} }{4d(\hat{H}_N)}\log\left(\frac{1}{N}\sum_{i=1}^{N} \frac{|X(\frac{i}{N}) - X(\frac{i-1}{N})|^2}{\hat{\sigma}_N^2 N^{-2\hat{H}_N}} \right) = m_NL_1(N)+ m_NL_2(N).
	\end{equation}
\autoref{thm: estimators convergence for sRP} guarantees that $m_N\xrightarrow[N\to\infty]{a.s.} 1$
and so the second summand in \eqref{2var non-convergence} converges to $Z^H(1)$ almost surely. However, \autoref{thm: estimators convergence for sRP} also implies that the first summand diverges because 
\begin{align*}
	N^{1-H} \log N |\hat{H}_N - H| \xrightarrow[N \to \infty]{a.s.} \infty,\\
	N^{1-H}|\log \hat{\sigma}_N - \log \sigma|\xrightarrow[N \to \infty]{a.s.} \infty.
\end{align*}
To make the second summand in \eqref{2var non-convergence} converge to zero, we can decelerate the growth factor $N^{1-H}$ in the 2-variation by taking only each $k$-th observation, where
\[k := \lfloor N/N^{\delta}\rfloor \quad \text{for some $0< \delta < 1$}, \quad (\lfloor.\rfloor\ \text{ is the floor function}),
\] 
when approximating $Z^{H}(1)$ while $H$ and $\sigma$ have to be estimated from all $N$ observations at the same time. Note that the increased time step becomes 
\begin{equation} 
  \label{eq: def h_N}
h_N = \frac{k}{N}  = \frac{\lfloor N/N^{\delta}\rfloor}{N} \sim \frac{1}{N^\delta}, 
\end{equation}
where we write $a_n\sim b_n$ whenever $\lim{(a_n/b_n)}=1$, and the number of observations available for approximating $Z^{H}(1)$ is
\[
n_N = \lfloor 1/h_N\rfloor \sim N^{\delta}.
\]
The following theorem confirms that the possibility to identify $Z^H(1)$ together with the speed of convergence, formulated in \autoref{Thm:V(Z) convergence} and \autoref{rem: speed of 2var to Z(1)}, remains valid also in this new situation where $n_N h_N$ may by less then 1. 

\begin{Theorem}
\label{Thm:decelerated V(Z) convergence}
Consider a sequence of positive numbers $\{h_N\}_{N=1}^{\infty}$ such that $n_N = \lfloor 1/h_N\rfloor \sim N^{\delta}$ for some $0< \delta< 1$. Then the convergence

\begin{equation} \label{eq:2var a.s. decelerated speed}
N^{\delta \alpha} \Bigg[\left(\frac{1}{h_N}\right)^{1-H} \frac{1}{4 d(H)} h_N \sum_{i=1}^{n_N}\Bigg(\frac{\left|Z^H(i h_N) - Z^H((i-1) h_N)\right|^2}{h_N^{2H}} - 1 \Bigg) - Z^H(1) \Bigg] \xrightarrow[N \to \infty]{a.s.} 0
\end{equation}
holds for 
\[ \alpha <  
\begin{cases}
H - \frac{1}{2}, & \quad H \in (\frac{1}{2},\frac{3}{4}],\\
1-H, & \quad H \in (\frac{3}{4},1).\\
\end{cases} \]
\end{Theorem}

\begin{proof}
Take rescaled $\tilde{h}_N = c_N h_N$ so that $\tilde{h}_N n_N = 1$. Clearly, the scaling factor $c_N = 1/(h_N n_N)$ tends to $1$ as $N \to \infty$. We then use $H$-self-similarity of the Rosenblatt process to obtain the equality of distributions 
\begin{equation} \label{eq: convergence of 2-var for h}
\begin{aligned}
&\left(\frac{1}{h_N}\right)^{1-H} \frac{1}{4 d(H)} h_N \sum_{i=1}^{n_N}\Bigg(\frac{\left|Z^H(i h_N) - Z^H((i-1) h_N)\right|^2}{h_N^{2H}} - 1 \Bigg) - Z^H(h_N n_N) \\
& \overset{d}{=} \frac{1}{c_N^H}\Bigg[\left(\frac{1}{\tilde{h}_N}\right)^{1-H} \frac{1}{4 d(H)} \tilde{h}_N \sum_{i=1}^{n_N}\Bigg(\frac{\left|Z^H(i \tilde{h}_N) - Z^H((i-1) \tilde{h}_N)\right|^2}{\tilde{h}_N^{2H}} - 1 \Bigg) - Z^H(1) \Bigg] \\
& = \frac{1}{c_N^H}\Bigg[ \frac{n_N^{1-H}}{4 d(H)} \frac{1}{n_N} \sum_{i=1}^{n_N}\Bigg(\frac{\left|Z^H(i/n_N) - Z^H((i-1)/n_N)\right|^2}{n_N^{-2H}} - 1 \Bigg) - Z^H(1) \Bigg] \\
& =\frac{1}{c_N^H}\Bigg[ \frac{n_N^{1-H}}{4d(H)}V_{n_N}(Z^H) - Z^H(1) \Bigg].
\end{aligned}
\end{equation}
Equation \eqref{eq:2var L2 speed} guarantees 
\[
\frac{n_N^{1-H}}{4d(H)}V_{n_N}(Z^H) - Z^H(1) = 
\begin{cases}
O_{L^2(\Omega)}(n_N^{\frac{1}{2}-H}), & \quad H \in (\frac{1}{2},\frac{3}{4}),\\
O_{L^2(\Omega)}(\sqrt{\log n_N} n_N^{\frac{1}{2}-H}), &\quad H = \frac{3}{4},\\
O_{L^2(\Omega)}(n_N^{H-1}), & \quad H \in (\frac{3}{4},1).\\
\end{cases}
\]
The equality of laws above implies the equality of the $L^2(\Omega)$-norms and so the same speed of $L^2(\Omega)$-convergence holds true for the first term in \eqref{eq: convergence of 2-var for h}. 
Stationarity of increments and $H$ self-similarity of the Rosenblatt process $Z^H$ further imply 
\[
\lVert  Z^H(h_N n_N) - Z^H(1) \rVert_{L^2(\Omega)} = (1- h_N n_N)^H \lVert Z^H(1) \rVert_{L^2(\Omega)} \leq (h_N )^H  \sim (n_N )^{-H} . 
\]
Altogether, we have 
\begin{equation} \label{eq: speed of L2 conv of 2-var for h}
\begin{aligned}
&\left(\frac{1}{h_N}\right)^{1-H} \frac{1}{4 d(H)} h_N \sum_{i=1}^{n_N}\Bigg(\frac{\left|Z^H(i h_N) - Z^H((i-1) h_N)\right|^2}{h_N^{2H}} - 1 \Bigg) - Z^H(1) \\
&  = 
\begin{cases}
O_{L^2(\Omega)}(N^{\delta(\frac{1}{2}-H)}), & \quad H \in (\frac{1}{2},\frac{3}{4}),\\
O_{L^2(\Omega)}(\sqrt{\log N} N^{\delta(\frac{1}{2}-H)}), &\quad H = \frac{3}{4},\\
O_{L^2(\Omega)}(N^{\delta(H-1)}), & \quad H \in (\frac{3}{4},1).\\
\end{cases}
\end{aligned}
\end{equation}
The claim can now be proved by same arguments as in the proof of \autoref{Thm:V(Z) convergence}.
\end{proof}
We can also rephrase the statements of \autoref{Lem: V(sZ) convergence} about convergence of the logarithmic version of the 2-variation with $N$ replaced by $1/h_N \sim N^{\delta}$. In comparison with the proof of \autoref{Lem: V(sZ) convergence}, an additional error term $h_N n_N - 1$ would appear inside the logarithm, but it is dominated by the original error term $r_N$. The rest of the proof would be analogous. Similarly, we can repeat the proof of \autoref{Thm:V(X) convergence} with slight modifications so that we obtain
\begin{equation}
\left(\frac{1}{h_N}\right)^{1-H} \frac{1}{4 d(H)} h_N \sum_{i=1}^{n_N}\left(\frac{\left|X\left(i h_N\right) - X\left((i-1) h_N\right)\right|^2}{\sigma^2 h_N^{2H}} - 1 \right) \quad \xrightarrow[N \to \infty]{a.s.}\quad Z^H(1),
\end{equation}
or, for logarithmic version,
\begin{equation}
\left(\frac{1}{h_N}\right)^{1-H} \frac{1}{4 d(H)} \log \left( h_N \sum_{i=1}^{n_N}\frac{\left|X\left(i h_N\right) - X\left((i-1) h_N\right)\right|^2}{\sigma^2 h_N^{2H}} \right) \quad \xrightarrow[N \to \infty]{a.s.}\quad Z^H(1).
\end{equation}

All these considerations motivate us to define the decelerated $2$-variation with estimated diffusion parameters by
	\[ \widehat{W}^{[0,1]}_{h_N}(X):= \left(\frac{1}{h_N}\right)^{1-\hat{H}_N}\frac{1}{4d(\hat{H}_N)} \log\left(h_N\sum_{i=1}^{n_N} \frac{\left|X(i h_N) - X((i-1) h_N)\right|^2}{\hat{\sigma}_N^2 h_N^{2\hat{H}_N}} \right), \]
that can be written as
	\begin{equation}
	\label{eq:dec_2var}
	\widehat{W}^{[0,1]}_{h_N}(X) = m_N^\delta L_1^\delta(N) + m_N^\delta L_2^\delta (N)
	\end{equation}
with 
	\begin{align*}
		m_N^\delta & := \frac{d(H)}{d(\hat{H}_N)} \e^{- (H-\hat{H}_N) \log h_N} \\
		L_1^\delta(N) & :=  \left(\frac{1}{h_N}\right)^{1-H}\frac{1}{4d(H)}  \left[ 2(\hat{H}_N - H) \log \frac{1}{h_N}  -2(\log \hat{\sigma}_N - \log \sigma) \right]\\
		L_2^\delta(N) & := \left(\frac{1}{h_N}\right)^{1-H} \frac{1}{4d(H)} \log \left(h_N \sum_{i=1}^{n_N}\frac{\left|X(i h_N) - X((i-1) h_N)\right|^2}{\sigma^2 h_N^{2H}}\right).
	\end{align*} 
Clearly, the second summand in \eqref{eq:dec_2var} converges almost surely to $Z^H(1)$ as $N\to \infty$, whereas for the first summand, we have:
\begin{equation}
\begin{aligned}
& \left(\frac{1}{h_N}\right)^{1-H} \left(\log \frac{1}{h_N}\right) (\hat{H}_N - H) \,\sim\, \delta N^{\delta (1-H)} (\log N) (\hat{H}_N - H)  \xrightarrow[N \to \infty]{a.s.} 0, \\
& \left(\frac{1}{h_N}\right)^{1-H} (\log \hat{\sigma}_N - \log \sigma) \,\sim\,N^{\delta (1-H)} (\log \hat{\sigma}_N - \log \sigma) \xrightarrow[N \to \infty]{a.s.} 0,
\end{aligned}
\end{equation}
by \autoref{thm: estimators convergence for sRP}. In result, the decelerated 2-variation with estimated diffusion parameters enjoys the desired convergence:
\[
\widehat{W}^{[0,1]}_{h_N}(X)  \xrightarrow[N \to \infty]{a.s.} Z^H(1).
\]
Now it follows immediately (by repeating the arguments from \autoref{subsec: Estimating drift}) that for
\[ \widehat{\mathbb{U}}_N^\delta :=
\begin{bmatrix} 
X(1/2) - X(0) - \hat{\sigma}_N \widehat{W}^{[0,1/2]}_{h_N}(X)\\
X(1) - X(1/2) - \hat{\sigma}_N \widehat{W}^{[1/2,1]}_{h_N}(X)
\end{bmatrix}
\]
we have
\[
\widehat{\mathbf{\varepsilon}}_N^\delta := \widehat{\mathbb{U}}_N^\delta - \lambda \cdot \mathbb{Z} \xrightarrow[N \to \infty]{a.s.} 
\begin{bmatrix} 
0 \\
0 \\
\end{bmatrix}.
\]
Consequently, we can define the plug-in estimator by
\begin{equation}\label{def: plug-in lambda} 
\hat{\lambda}_{N}^\delta := (\mathbb{Z}_{h_N}^T \cdot \mathbb{Z}_{h_N})^{-1} \cdot \mathbb{Z}_{h_N}^T \cdot \widehat{\mathbb{U}}_N^\delta,
\end{equation}
where $\mathbb{Z}_{h_N}$  is the discrete approximation (by using Riemann sums) of $\mathbb{Z}$ with time step $h_N$. Applying the same arguments as in \autoref{subsec: Estimating drift} provides us with strong consistency of this estimator for any $0<\delta<1$:
\begin{equation}\label{eq: plug-in lambda consistent} 
\hat{\lambda}_{N}^\delta \xrightarrow[N \to \infty]{a.s.} \lambda.
\end{equation}
It is natural to ask which $\delta \in (0,1)$ should one choose. In what follows, we will find such $\delta$ that optimizes the speed of convergence in \eqref{eq: plug-in lambda consistent}. Recall (see \autoref{subsec: Estimating drift} for details) that the speed of convergence of the drift estimator is determined by the speed of convergence of the 2-variation towards $Z^{H}(1)$. Therefore, we will examine the speed of this convergence for the decelerated 2-variation for different values of $\delta$. In particular, we wish to maximize $a$ for which
\begin{equation}\label{plug-in required rate}
N^{a}\left| \widehat{W}^{[0,1]}_{h_N}(X) - Z^H(1) \right| \xrightarrow[N \to \infty]{a.s.} 0.
\end{equation}
Initially, we can estimate the left-hand side of \eqref{plug-in required rate} similarly to \eqref{eq:dec_2var} %\eqref{eq: dec 2var convergence} 
as follows:
\[ N^{a}\left| \widehat{W}^{[0,1]}_{h_N}(X) - Z^H(1) \right| \leq 
\underbrace{N^a |m_N^\delta L_1^\delta(N)|}_{=:S_1} + \underbrace{N^a |m_N^\delta L_2^\delta(N) - Z^H(1)|}_{=:S_2}.\]
The sufficient (and necessary) condition for the first summand $S_1$ to converge to zero follows directly from \eqref{eq: convergence H sigma} and it reads:
\begin{equation}\label{convergence cond1}
a + \delta (1-H) < 1-H, \quad \text{or, equivalently, }\quad a <(1- \delta) (1-H). 
\end{equation}
For the second summand $S_2$, we can write
	\begin{equation*}
		S_2  \leq \underbrace{N^a |L_2^\delta(N) - Z^H(1)|}_{=:S_{2,1}} + \underbrace{N^a |m_N^\delta-1| |L_2^\delta(N)|}_{=:S_{2,2}}.
	\end{equation*}
Now, a sufficient condition for the convergence of $S_{2,1}$ to zero follows directly from a version of convergence \eqref{eq:2var a.s. decelerated speed} for the logarithmic 2-variation  and it reads as follows:
\begin{equation}\label{convergence cond2}
a = \delta \alpha <
\begin{cases}
\delta(H - \frac{1}{2}), & \quad H \in (\frac{1}{2},\frac{3}{4}],\\
\delta(1-H), & \quad H \in (\frac{3}{4},1).\\
\end{cases}
\end{equation}
For the convergence of $S_{2,2}$ to zero, write
\[
N^{a} \left| m_N^\delta  - 1 \right| = N^{a} \left| \e^{x_N}  - 1 \right|
\]
with 
	\[x_N:=[\log(d(H)) -\log(d(\hat{H}_N))] + (H - \hat{H}_N) \log \frac{1}{h_N}\]
so that we see that \[N^\alpha|m_N^\delta-1| = N^{\alpha} \left| x_N + O(x_N^2) \right|, \quad N\to\infty.\]
Since the function $\log(d(\cdot))$ is continuously differentiable with positive (and finite) derivative in the interval $(1/2, 1)$, and $1/h_N \sim N^\delta$, we have
\[
x_N = O(|H - \hat{H}_N| \log N), \quad N\to\infty.
\]
The convergence rate for $\hat{H}_N$ in \eqref{eq: convergence H sigma} ensures the almost sure the convergence of $S_{2,2}$ to zero as long as $\alpha < 1-H$. This last condition is, however, guaranteed by condition \eqref{convergence cond1}.
Thus, we obtain the required convergence in \eqref{plug-in required rate} if the two conditions \eqref{convergence cond1} and \eqref{convergence cond2} are satisfied.

\begin{Remark}
Observe that the more intensive the deceleration is (i.e. the smaller $\delta$ we take), the faster the convergence of term $S_1$ to zero is obtained (in fact, the divergence of this term in \eqref{2var non-convergence} was the reason for the introduction of the deceleration), but the worse convergence rate for term $S_2$ (which expresses convergence of the decelerated 2-variation towards $Z^H(1)$). This behavior is therefore in line with what is expected and it enables us to find the optimal $\delta$ which would balance the two convergence rates.
\end{Remark}

Given $H \in (1/2 , 1)$, we now wish to find such $\delta \in (0,1)$ that maximizes $a$ for which conditions \eqref{convergence cond1} and \eqref{convergence cond2} hold. It follows that the optimum deceleration rate is
\begin{equation}\label{eq: delta^o}
\delta^o = 
\begin{cases}
2(1-H), & \quad H \in (\frac{1}{2},\frac{3}{4}], \\
\frac{1}{2}, & \quad H \in (\frac{3}{4},1), \\
\end{cases}
\end{equation}
and the corresponding convergence rate in \eqref{plug-in required rate} is
\begin{equation}\label{eq: alpha^o}
a^o < 
\begin{cases}
2(1-H)(H - \frac{1}{2}), & \quad H \in (\frac{1}{2},\frac{3}{4}], \\
\frac{1}{2}(1-H), & \quad H \in (\frac{3}{4},1).
\end{cases}
\end{equation}

The results on the speed of convergence of the optimized plug-in estimator can now be given.

\begin{Theorem}\label{thm: plug-in lambda optim convergence}
The optimized plug-in estimator 
\begin{equation}\label{def: optim lambda} 
\hat{\lambda}_{N}^{\delta^o} = (\mathbb{Z}_{h_N}^T \cdot \mathbb{Z}_{h_N})^{-1} \cdot \mathbb{Z}_{h_N}^T \cdot \widehat{\mathbb{U}}_N^{\delta^o},
\end{equation}
    with $\delta^o$ defined in \eqref{eq: delta^o} and the corresponding $h_N$ in \eqref{eq: def h_N}, is strongly consistent, i.e.
\begin{equation*}\label{eq: plug-in lambda optim consistent}
\hat{\lambda}_{N}^{\delta^o} \xrightarrow[\substack{N\in 2\mathbb{N}\\ N\to\infty}]{a.s.} \lambda.
\end{equation*}
The corresponding convergence rate is at least
\begin{equation*}\label{eq: lug-in lambda optim convergence}
N^{a^o}\left|\hat{\lambda}_{N}^{\delta^o}(\omega) - \lambda \right| \xrightarrow[\substack{N\in 2\mathbb{N}\\ N\to\infty}]{a.s.} 0,
\end{equation*}
with $a^o$ specified in \eqref{eq: alpha^o}.
\end{Theorem}

\begin{proof}
In order to repeat the arguments from the proof of \autoref{thm: hat(lambda)1 convergence}, recall that for \[\widehat{\mathbf{\varepsilon}}_N^{\delta^o} = \widehat{\mathbb{U}}_N^{\delta^o} - \lambda \cdot \mathbb{Z}\] we have 
\begin{align*}
N^{a^o}\widehat{\mathbf{\varepsilon}}_N^{\delta^o} & = N^{a^o}(\widehat{\mathbb{U}}_N^{\delta^o} - \lambda \cdot \mathbb{Z}) \\
&= N^{a^o} \begin{bmatrix} 
\sigma Z^H(1/2) - \hat{\sigma}_N \widehat{W}^{[0,1/2]}_{h_N}(X)\\
\sigma (Z^H(1) - Z^H(1/2)) - \hat{\sigma}_N \widehat{W}^{[1/2,1]}_{h_N}(X)
\end{bmatrix} \\
&= N^{a^o} \begin{bmatrix} 
\sigma \left( Z^H(1/2) - \widehat{W}^{[0,1/2]}_{h_N}(X)\right) - \left(\hat{\sigma}_N - \sigma \right) \widehat{W}^{[0,1/2]}_{h_N}(X)\\
\sigma \left((Z^H(1) - Z^H(1/2)) - \widehat{W}^{[1/2,1]}_{h_N}(X)\right) - \left(\hat{\sigma}_N - \sigma \right) \widehat{W}^{[1/2,1]}_{h_N}(X)
\end{bmatrix}.
\end{align*}
It follows by \eqref{plug-in required rate} and \eqref{eq: convergence H sigma} that
\[
N^{a^o}\widehat{\mathbf{\varepsilon}}_N^{\delta^o}  \xrightarrow[\substack{N\in 2\mathbb{N}\\ N\to\infty}]{a.s.} 
\begin{bmatrix} 
0 \\
0 \\
\end{bmatrix}.
\]
Further note that convergence rate of $\mathbb{Z}_{h_N}$ analogical to \eqref{eq:Z_N convergence rate} guarantees that for almost every $\omega$, the term
\begin{equation*}
N^{a^o} \left| \mathbb{Z}_{h_N}(\omega) - \mathbb{Z}(\omega) \right| \leq C(\omega) (h_N)^{\gamma} N^{a^o} \sim \frac{N^{a^o}}{N^{\delta^o \gamma}}
\end{equation*}
converges to zero for all $\gamma<H$. Now the proof is finished by repeating the arguments from the proof of \autoref{thm: hat(lambda)1 convergence}.
\end{proof}

\begin{figure} 
    \includegraphics[scale=0.35, keepaspectratio]{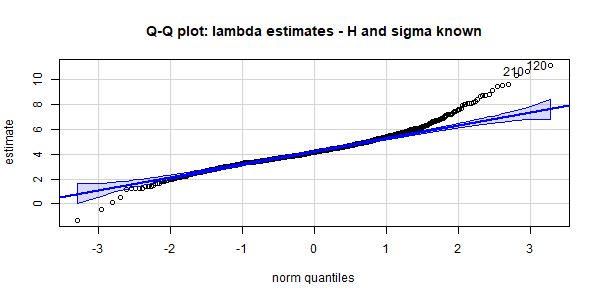}
    \includegraphics[scale=0.35, keepaspectratio]{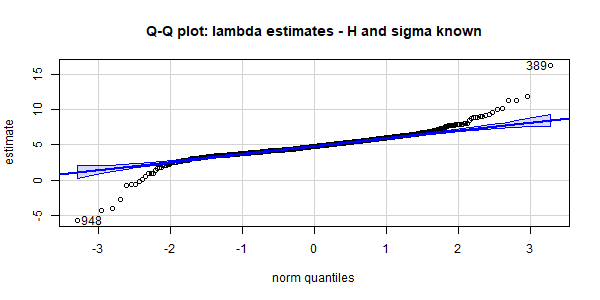}
    \includegraphics[scale=0.35, keepaspectratio]{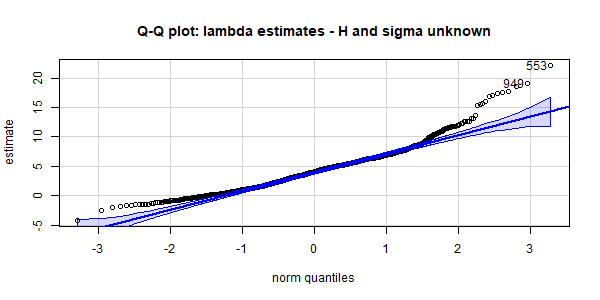}
    \includegraphics[scale=0.35, keepaspectratio]{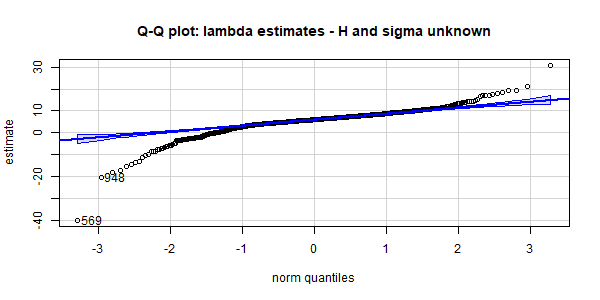}
\caption{Q-Q plots comparing sample quantiles of drift estimates $\hat{\lambda}_{1,N}$ (top) and $\hat{\lambda}_{N}^{\delta^o}$ (bottom) with Gaussian quantiles (incl. 95\% confidence bands) for equations (ROU) on the left and (RSDE) on the right.}
  \label{fig:L QQplot}
\end{figure}

\begin{Remark}
Both estimators $\hat{\lambda}_{1,N}$ (assumes $H$ and $\sigma$ known) and the decelerated plug-in estimator $\hat{\lambda}_{N}^{\delta^o}$ (uses estimates $\hat{H}_N$ and $\hat{\sigma}_N$ of $H$ and $\sigma$, respectively, as well as the decelerated 2-variation) are strongly consistent. However, the convergence rate of the former is better than that of the latter. 

Theoretical derivation of asymptotic distributions of $\hat{\lambda}_{1,N}$  and $\hat{\lambda}_{N}^{\delta^o}$ are beyond the scope of this paper and may be subject of further research. The simulation results indicate non-Gaussianity of the distribution with heavy tails, especially the upper tail. For details, see the comparison of sample quantiles of the (simulated) distributions of $\hat{\lambda}_{1,N}$ and $\hat{\lambda}_{N}^{\delta^o}$ with Gaussian quantiles in Figure \ref{fig:L QQplot}. Similar Monte Carlo or bootstrap techniques might be also used to construct approximate confidence intervals in our situation when the true (asymptotic) distributions of the estimators of $\lambda$ are not known but it is possible to sample from the distribution.     

\end{Remark}

\begin{Remark} Let us briefly comment on the case of general Hermite processes. Denote by $\{Z^{q,H}(t)\}_{t \in [0,1]}$ an Hermite process of order $q \in \mathbb{N}$ with Hurst parameter $H \in (1/2,1)$ (see e.g. \cite{ChronopoulouEtAl11}, Definition 2.1 for a rigorous definition). Since this process is defined as the Wiener-It\^o multiple integral of order $q$ with respect to the standard Wiener process, random variables $Z^{q,H}(t)$, $t\in\R$, belong to the $q$-th Wiener chaos. Moreover, the process is centered with stationary increments, self-similar with index $H$, it has $\gamma$-H\"older continuous trajectories for every $\gamma<H$, and it possesses the same covariance function as the fractional Brownian motion, i.e.
\[
\mathbb{E}\left[Z^{q,H}(s) Z^{q,H}(t)\right] = \frac{1}{2}(t^{2H} + s^{2H} - |t-s|^{2H}), \quad s,t\in [0,1].
\]  
The class of Hermite processes represents a natural generalization to that of fractional Brownian motions (which is an Hermite process with $q=1$) and to the Rosenblatt process (which is an Hermite process with $q=2$).
A detailed study of the asymptotic properties of 2-variations of Hermite processes $V_N(Z^{q,H})$ (see \eqref{def:2variation} for the precise definition of $V_N(\cdot)$) can be found in the article \cite{ChronopoulouEtAl11}. Let us briefly summarize that the authors decomposed the 2-variation into the individual Wiener chaoses as follows: 
\[
V_N(Z^{q,H}) = T_{2q} + c_{2q-2} T_{2q-2} + \ldots c_{4} T_{4} + c_{2} T_{2}, \quad \text{where } c_{2q-2k}=k! \binom{q}{k}^2.
\]
They found out that the $2$\textsuperscript{nd} chaos component $T_2$ is dominant in the decomposition and that there is the convergence
\begin{equation}
	\label{eq:L2 convergence of V_N for Hermit}
\lim_{N \to \infty} N^{2-2H'} V_N(Z^{q,H}) = \lim_{N \to \infty} N^{2-2H'} c_{2} T_{2} =  c_{1,H}^{1/2} c_2 Z^{2,2H'-1}(1),
\end{equation}
understood in the $L^2(\Omega)$ sense. Here, $H' = 1 + (H-1)/q$, $c_{1,H}$ is an explicit constant and $Z^{2,2H'-1}(1)$ is the standard Rosenblatt random variable corresponding to the Rosenblatt process with Hurst parameter $2H'-1$ evaluated at $t=1$. This Rosenblatt process is constructed from the same Wiener process as the original Hermite process $Z^{q,H}$. In addition, the individual components exhibit the following convergence:
\begin{equation}
	\label{eq: convergence of 2-variation components for Hermit}
	N^{(2-2H')(q-k)}T_{2q-2k} \quad\xrightarrow[N \to \infty]{L^2(\Omega)}\quad z_{k,h} Z^{2q-2k,\tilde{H}}, \quad \text{for every } k=1,\ldots, q-2,
\end{equation}
where  $\tilde{H} = (2q-2k)(H'-1) + 1$, $z_{k,h}$ is an explicit constant and $Z^{2q-2k,\tilde{H}}$ denotes an Hermite random variable of order $2q-2k$ with self-similarity index $\tilde{H}$.
The convergence \eqref{eq:L2 convergence of V_N for Hermit} can be utilized to identify parameters $\sigma$ and $H$ for a scaled Hermite process. On the other hand, the approach from section \ref{sec:scaled RP} would have to be adapted from the almost sure to $L^2(\Omega)$ convergence and this is outside the scope of the present work.\\
As far as the identification of $\sigma$ and $H$ for SDEs of type \eqref{eq:non-lin_SDE} driven by an Hermite process is concerned, it seems likely that the procedure outlined above would also be possible.
On the other hand, the procedure for estimating drift described above cannot be applied to SDEs driven by Hermite processes of an order higher than 2. This is because in such a case, the evaluation of the 2-variation of the trajectory of the solution does not identify the value of the corresponding noise that lives in the $q$-th chaos but only its component in the $2$\textsuperscript{nd} chaos. This component represents a Rosenblatt process (with a potentially different Hurst parameter) that is measurable with respect to the same Wiener process that generates the Hermite noise. Hence, the problem of singularity or equivalence of the distributions of the solutions to Hermite-driven SDEs with different drifts remains, to the best of our knowledge, open. 
\end{Remark}

\section{Simulations}\label{sec:Simulations}
In this section, we illustrate the theoretical results on asymptotic behavior of the estimators on simulated data. Consider the general equation
\begin{equation}
	\label{eq:non-lin_SDE simulations}
	 X(t) = X_0 + \lambda \int_0^t f(X(s))\d{s} + Y(t), \qquad t\in [0,1],
\end{equation}
with the following three specifications:
\begin{enumerate}[wide=0pt, leftmargin=*]
    \item[(ROU)] Rosenblatt Ornstein-Uhlenbeck process: \\
    $X_0 = 0.5, \lambda = 5, f(x) = -x$, and $Y = Z^H$ (Rosenblatt process), $H = 0.75$,\\
    \item[(RSDE)] Non-linear SDE driven by Rosenblatt process:\\
    $X_0 = 0, \lambda = 5, f(x) = x (1-x) (1+x)$, and $Y = Z^H$ (Rosenblatt process), $H = 0.75$,\\
    \item[(FOU)] Fractional Ornstein-Uhlenbeck process: \\
    $X_0 = 0.5, \lambda = 5, f(x) = -x$, and $Y = B^H$ (fractional Brownian motion), $H = 0.75$.\\
\end{enumerate}
The FOU process was chosen to investigate the effect of noise misspecification on the estimators. Simulations of the solutions to the above SDEs were performed in Wolfram Mathematica 13.1.\ via the Euler–Maruyama method (mesh size $\Delta t = 1/51200$). The underlying FBM was generated by the internal procedure  ``FractionalBrownianMotionProcess'' and the Rosenblatt process was simulated according to the method described in \cite{BarTud10} (originally introduced in \cite{AbrPip06}). The evaluation of the estimators and the subsequent analysis of the results was performed in the statistical software~R.

\begin{figure}
   \includegraphics[scale=0.17, keepaspectratio]{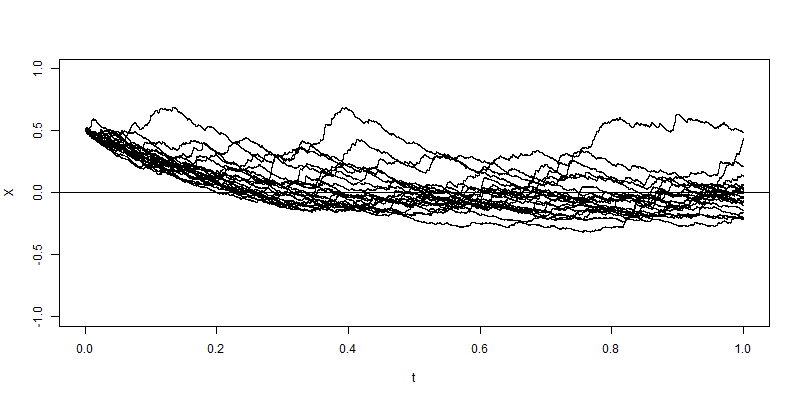} 
   \includegraphics[scale=0.17, keepaspectratio]{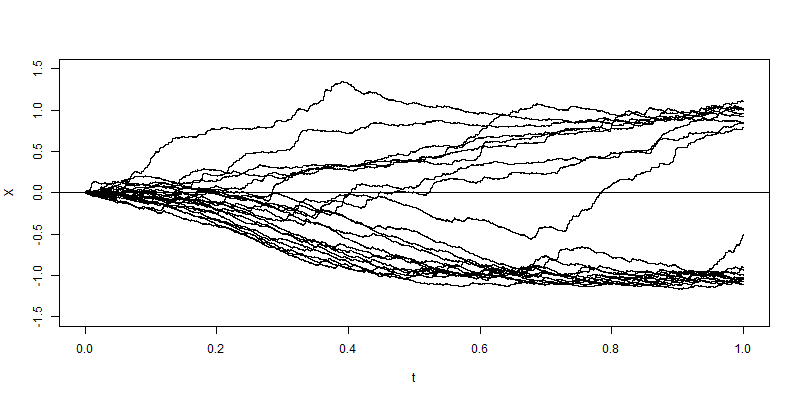}
   \includegraphics[scale=0.17, keepaspectratio]{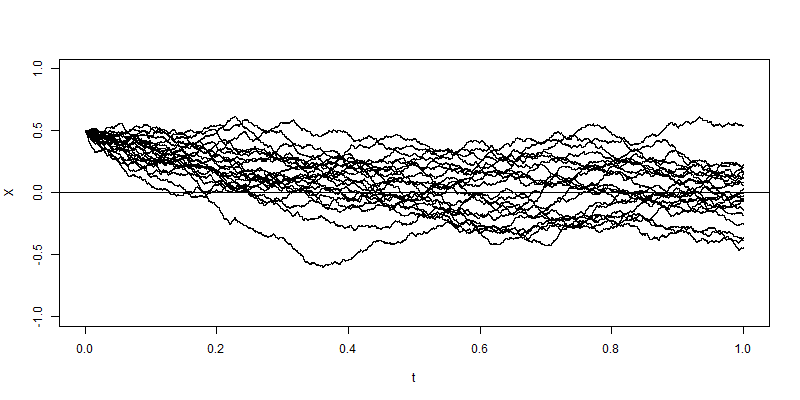}
  \caption{Sample solutions to the equations (ROU), (RSDE) and (FOU) ordered from left to right.}
  \label{fig:sample_solutions}
\end{figure}

We generated $1000$ trajectories for each of the three models above. Selected sample trajectories are shown in Figure \ref{fig:sample_solutions}. Subsequently, we evaluated the four estimators studied above for each trajectory. In particular, we calculated the Hurst parameter and the noise intensity estimates $\hat{H}_N$ and $\hat{\sigma}_N$, respectively, defined in \autoref{thm: estimators convergence for sRP}, as well as the drift parameter estimates $\hat{\lambda}_{1,N}$, defined by formula \eqref{def: hat(lambda)1} assuming $H$ and $\sigma$ known, and, finally, the optimized plug-in estimator $\hat{\lambda}_{N}^{\delta^o}$ defined in \eqref{def: optim lambda} without the knowledge of $H$ and $\sigma$. The sample values of the estimates were then summarized in boxplots and log-log plots of root-mean-square-errors (RMSE) for different values of $N$.

\begin{figure} 
    \includegraphics[scale=0.23, keepaspectratio]{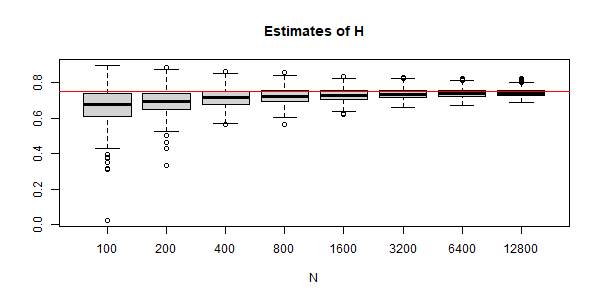}
    \includegraphics[scale=0.23, keepaspectratio]{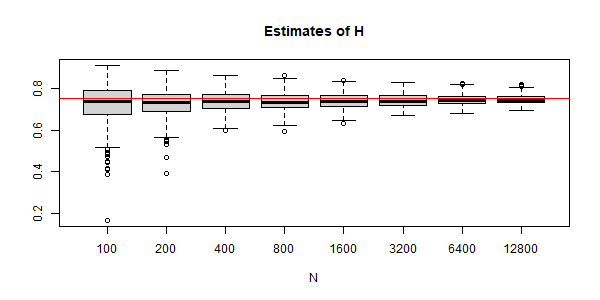}
    \includegraphics[scale=0.23, keepaspectratio]{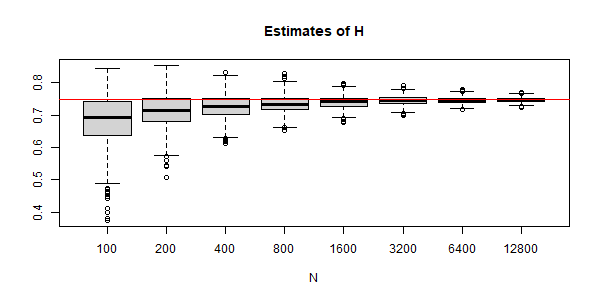}
    \includegraphics[scale=0.23, keepaspectratio]{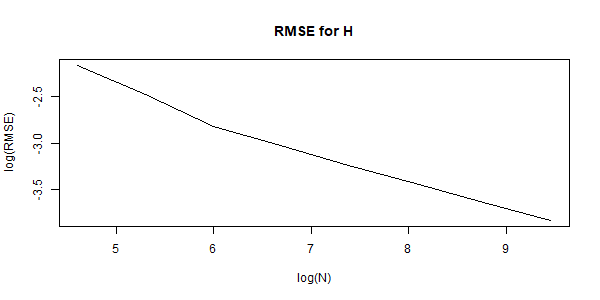}
    \includegraphics[scale=0.23, keepaspectratio]{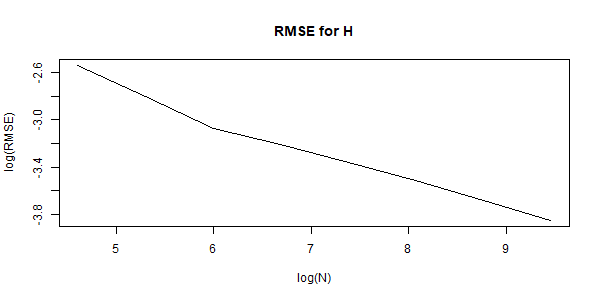}
    \includegraphics[scale=0.23, keepaspectratio]{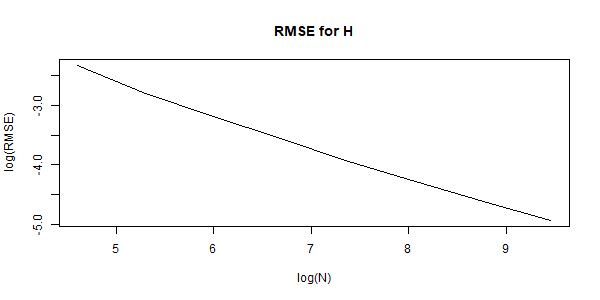}
\caption{Boxplots of estimates $\hat{H}_N$ of Hurst parameter $H$ (top) and log-log plots of corresponding RMSE (bottom) for equations (ROU), (RSDE) and (FOU) (ordered from left to right). Red lines represent true parameter values.}
  \label{fig:H estimates}
\end{figure}

\begin{figure} 
    \includegraphics[scale=0.23, keepaspectratio]{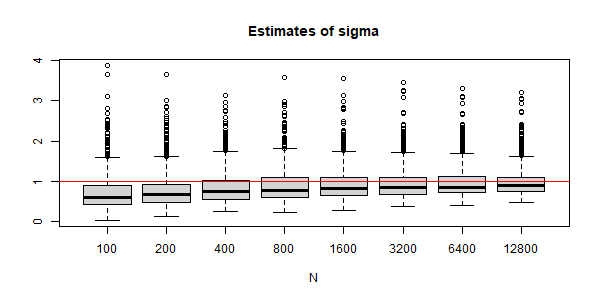}
    \includegraphics[scale=0.23, keepaspectratio]{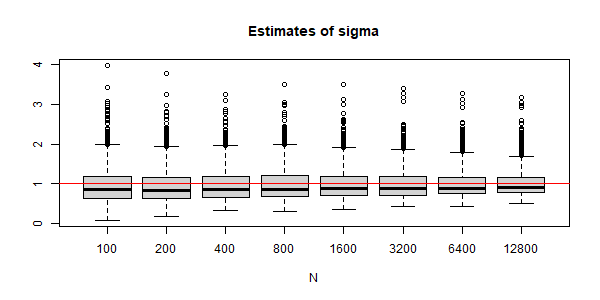}
    \includegraphics[scale=0.23, keepaspectratio]{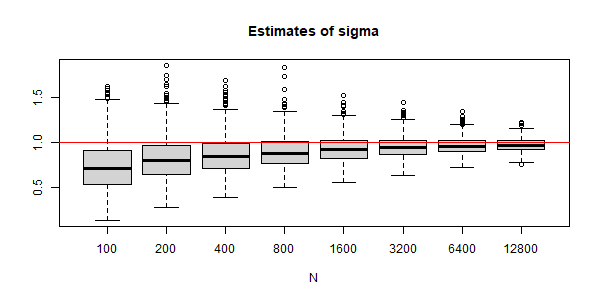}
    \includegraphics[scale=0.23, keepaspectratio]{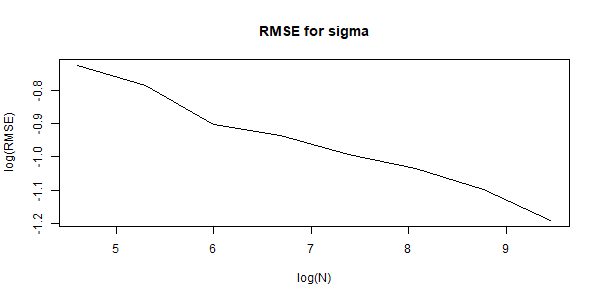}
    \includegraphics[scale=0.23, keepaspectratio]{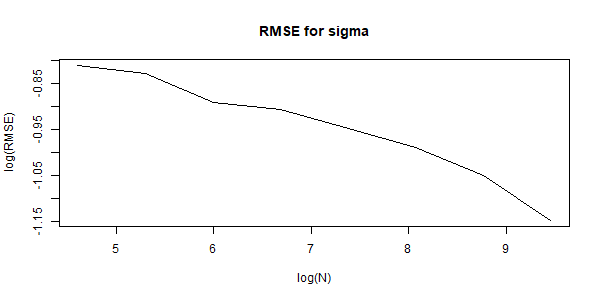}
    \includegraphics[scale=0.23, keepaspectratio]{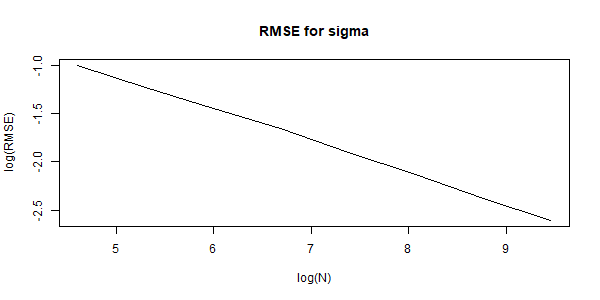}

\caption{Boxplots of estimates $\hat{\sigma}_N$ of noise intensity $\sigma$ (top) and log-log plots of corresponding RMSE (bottom) for equations (ROU), (RSDE) and (FOU) (ordered from left to right). Red lines represent true parameter values.}
  \label{fig:S estimates}
\end{figure}

The behavior of the estimator $\hat{H}_N$ is illustrated in Figure \ref{fig:H estimates}. We see clear convergence towards the true value for all three models. The limiting sample distributions for models (ROU) and (RSDE), driven by Rosenblatt process, are right-skewed. This corresponds to the Rosenblatt distribution claimed in \autoref{thm: estimators convergence for sRP}. The estimator converges also for the Ornstein-Uhlenbeck process driven by a FBM (FOU) but with a different (symmetric) limiting distribution. This robustness appears because the construction of the estimator builds mainly upon self-similarity of the noise, while the specific properties of the Rosenblatt process, summarized in \autoref{Thm:V(Z) convergence}, are only used to determine the asymptotic distribution. The same conclusions apply also for the estimator $\hat{\sigma}_N$ whose performance can be found in Figure \ref{fig:S estimates}. Somewhat slower convergence of $\hat{\sigma}_N$ compared to $\hat{H}_N$ may partly be explained by the additional factor $1/\log N$ in \autoref{thm: estimators convergence for sRP} and partly by possible imperfections in RP simulations (compare convergence in the models driven by RP with the model driven by FBM which is much easier to simulate). 

\begin{figure} 
    \includegraphics[scale=0.23, keepaspectratio]{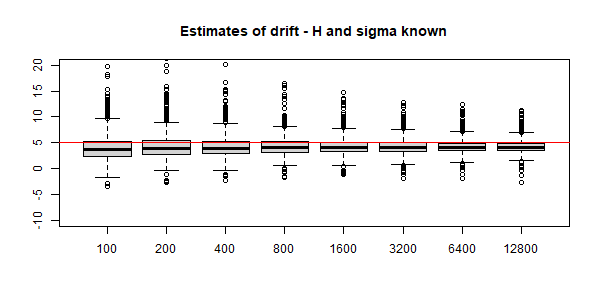}
    \includegraphics[scale=0.23, keepaspectratio]{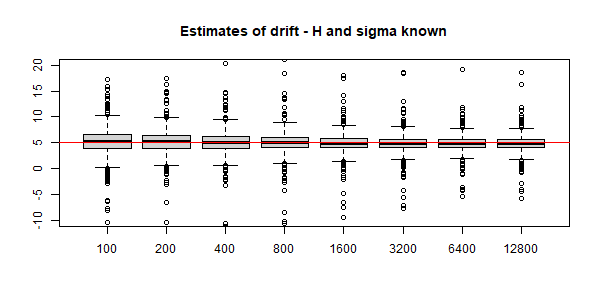}
    \includegraphics[scale=0.23, keepaspectratio]{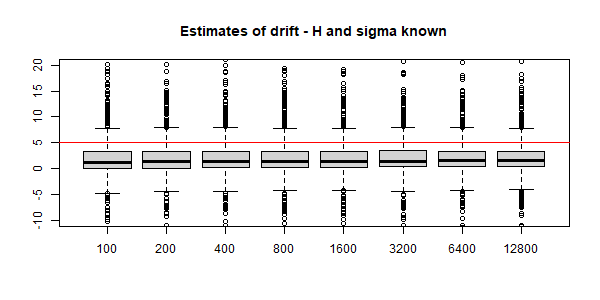}
    \includegraphics[scale=0.23, keepaspectratio]{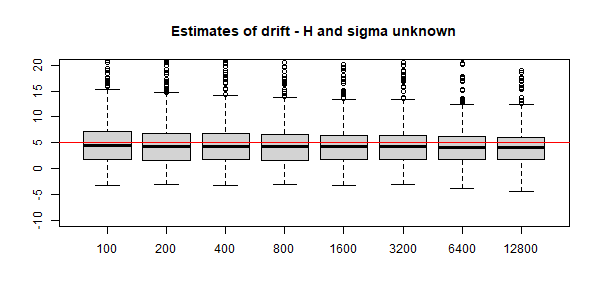}
    \includegraphics[scale=0.23, keepaspectratio]{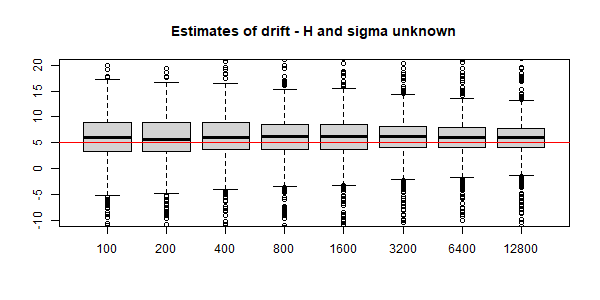}
    \includegraphics[scale=0.23, keepaspectratio]{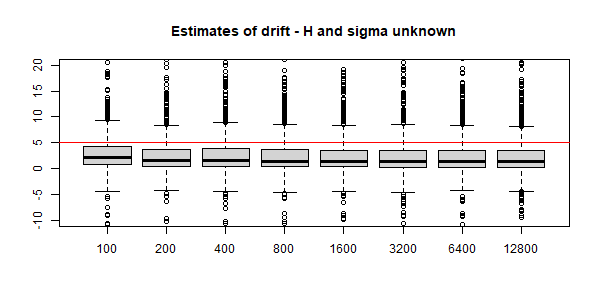}
\caption{Boxplots of drift estimates  $\hat{\lambda}_{1,N}$  considering $H$ and $\sigma$ known (top), and the optimized plug-in estimator $\hat{\lambda}_{N}^{\delta^o}$ with $H$ and $\sigma$ unknown (bottom) for equations (ROU), (RSDE) and (FOU) ordered from left to right. Red line represents true values.}
  \label{fig:L boxplots}
\end{figure}

\begin{figure} 
    \includegraphics[scale=0.23, keepaspectratio]{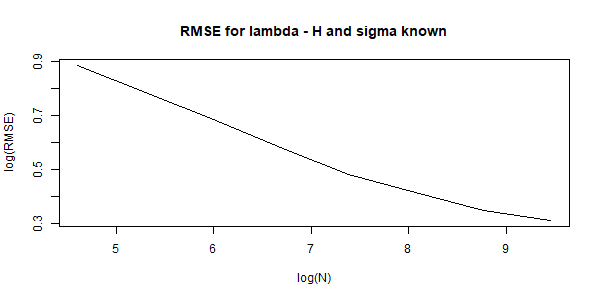}
    \includegraphics[scale=0.23, keepaspectratio]{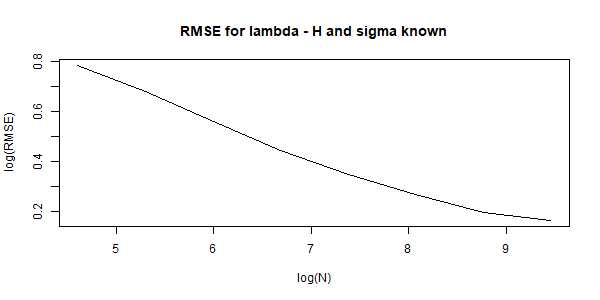}
    \includegraphics[scale=0.23, keepaspectratio]{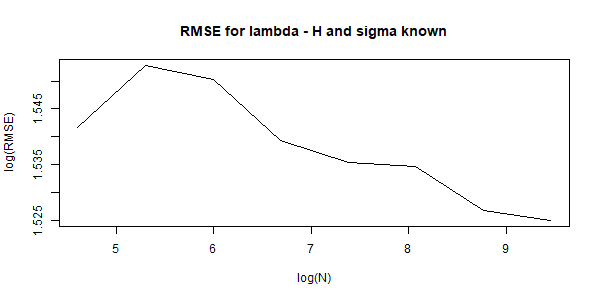}
    \includegraphics[scale=0.23, keepaspectratio]{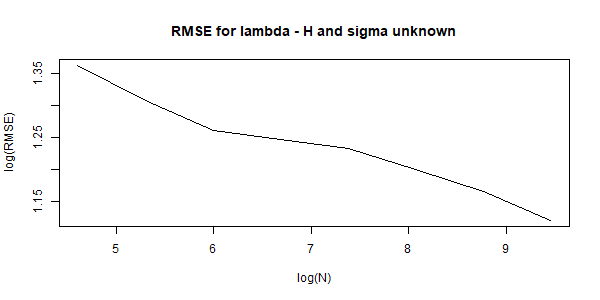}
    \includegraphics[scale=0.23, keepaspectratio]{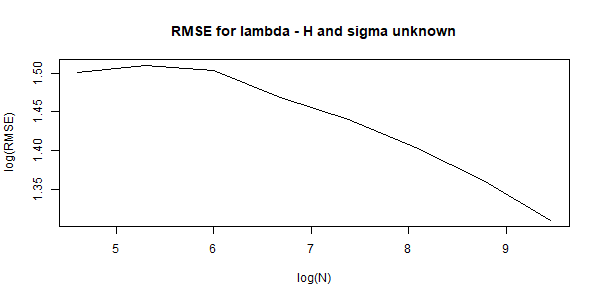}
    \includegraphics[scale=0.23, keepaspectratio]{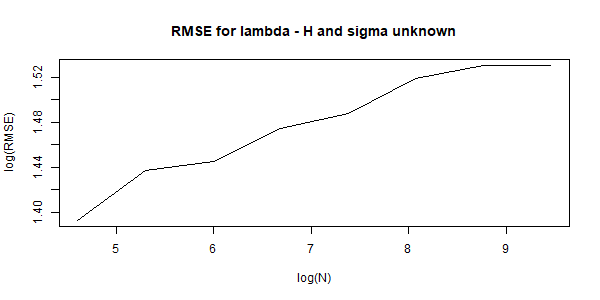}
\caption{Log-log plots of sample RMSE for drift estimates $\hat{\lambda}_{1,N}$ (top) and $\hat{\lambda}_{N}^{\delta^o}$ (bottom) for equations (ROU), (RSDE) and (FOU) ordered from left to right.}
  \label{fig:L RMSE}
\end{figure}

Both estimators of the drift parameter $\lambda$ are shown in Figures \ref{fig:L boxplots} and \ref{fig:L RMSE}. Correspondingly to \autoref{thm: hat(lambda)1 convergence} and \autoref{thm: plug-in lambda optim convergence}, we observe the convergence of both estimators for all models driven by a Rosenblatt process -- (ROU) and (RSDE), with $\hat{\lambda}_{1,N}$ converging faster then the plug-in estimator $\hat{\lambda}_{N}^{\delta^o}$. The observed speed of convergence is, however, lower compared to the theoretical one, which may probably be attributed to imperfections in RP simulations again. On the contrary, the two estimators fail to identify the drift parameter in the case when the model is driven by an FBM -- (FOU). This striking contrast illustrates the fundamental difference between SDEs driven by a Rosenblatt process and those driven by an FBM when it comes to drift identification. Whereas the former enables consistent estimation on a finite time interval (due to singularity of measures), consistent estimation in such setting is impossible for the latter (due to equivalence of measures).

\bibliographystyle{plain}
\bibliography{bibi}

\end{document}